\documentclass[11pt,a4paper,reqno]{amsart}
\ifdefined\pdfmapfile\pdfmapfile{+lm.map}\pdfmapfile{+euler.map}\pdfmapfile{+symbols.map}\pdfmapfile{+cm.map}\pdfmapfile{+cmextra.map}\fi
\usepackage[T1]{fontenc}
\usepackage{lmodern,amsmath,amssymb,mathtools,booktabs}
\usepackage[margin=27mm]{geometry}
\usepackage{hyperref}
\hypersetup{pdfborder={0 0 0},pdftitle={Approximation properties of Banach space pairs and quotients of lp}}
\newtheorem{theorem}{Theorem}[section]
\newtheorem{lemma}[theorem]{Lemma}
\newtheorem{proposition}[theorem]{Proposition}
\newtheorem{corollary}[theorem]{Corollary}
\theoremstyle{definition}
\newtheorem{problem}{Problem}
\theoremstyle{remark}

\numberwithin{equation}{section}
\newcommand{\E}{\mathbb E}
\newcommand{\tr}{\operatorname{tr}}

\newcommand{\1}{\mathbf1}

\newcommand{\R}{\mathbb R}
\newcommand{\ns}[2]{|#1|_{[#2]}}
\newcommand{\eu}[1]{|#1|_2}
\title[Approximation properties of pairs and quotients]{Approximation properties of Banach space pairs and quotients of $\ell_p$}

\author{Gang Cai$^\dag$}
\address[G. Cai]{School of Mathematical Sciences, Chongqing Normal University, Chongqing 401331, China}
\email{caigang-aaaa@163.com}

\author[Longfa Sun]{Longfa Sun}
\address[L. Sun]{Hebei Key Laboratory of Physics and Energy Technology, School of Mathematics and Physics, North China Electric Power University, Baoding 071003, China}
\email{sun.longfa@ncepu.edu.cn}

\thanks{$^\dag$ Corresponding author.\\
This work is supported by the National Natural Science Foundation of China (Grant No. 12571141), the
Chongqing Talent Support Program (Grant No. cstc2024ycjh-bgzxm0121), and the Fundamental Research Funds for the Central Universities (Grant No. 2025MS178).}

\date{}
\ifcsname subjclassname@2020\endcsname\else
  \expandafter\def\csname subjclassname@2020\endcsname{\textup{2020} Mathematics Subject Classification}
\fi
\begin{document}
\begin{abstract}
We construct a uniformly convex Banach space $X$ with a $1$-unconditional basis and a closed subspace $Y$ such that each of $X$, $Y$, and $X/Y$ has the metric approximation property, yet the pair $(X,Y)$ fails the approximation property (AP). Both $Y$ and $X/Y$ admit finite-dimensional decompositions. For every $2<p<\infty$, we construct an isomorphic embedding $J_p:\ell_p\to\ell_p$ such that $Q_p=\ell_p/J_p(\ell_p)$ fails the AP. For the conjugate exponent $q$, the adjoint $J_p^*:\ell_q\to\ell_q$ is surjective and its kernel also fails the AP. For each fixed $p$, every infinite-dimensional $\mathcal L_p$-space $V$ contains a closed subspace $W$ with the uniform approximation property such that $W\simeq V$ and $V/W\simeq Q_p$. These constructions give negative answers to the questions of Figiel, Johnson and Pe{\l}czy\'nski on approximation for pairs and quotients, including Szankowski's question on the bounded approximation property of quotients of $\ell_p$.
\end{abstract}
\maketitle
\section{Introduction}
Throughout this paper, all Banach spaces are real or complex. Write $\mathcal F(X)$ for the space of bounded finite-rank operators on $X$. A Banach space $X$ has the approximation property (AP) if, for every compact $K\subset X$ and every $\varepsilon>0$, there is a $T\in\mathcal F(X)$ such that
\[
                  \sup_{x\in K}\|Tx-x\|<\varepsilon.
\]
If $T$ can always be chosen with $\|T\|\le\lambda$, then $X$ has the $\lambda$-bounded approximation property ($\lambda$-BAP). The bounded approximation property (BAP) means the $\lambda$-BAP for some $\lambda<\infty$, and the metric approximation property (MAP) is the $1$-BAP. The uniform approximation property (UAP) requires a constant $C$ and a function $r:\mathbb N\to\mathbb N$ such that for every $k$-dimensional subspace $F\subset X$ there is a $T\in\mathcal F(X)$ satisfying
\[
       T|_F=I_F,\qquad \|T\|\le C,\qquad \operatorname{rank}T\le r(k).
\]
In particular the bound on the rank is independent of the position of $F$ in $X$. Both the MAP and the UAP imply the BAP, and the BAP implies the AP.

A linear map $J:X\to Z$ is an \emph{isomorphic embedding} if there are constants $0<c\le C<\infty$ such that
\[
               c\|x\|_X\le\|Jx\|_Z\le C\|x\|_X\qquad (x\in X).
\]
The lower bound implies that $J$ is injective and its range is closed, since $X$ is complete. Thus $J:X\to J(X)$ has a bounded inverse. We call $J(X)$ a copy of $X$ in $Z$. If $J(X)=Z$, we write $X\simeq Z$. The AP, BAP and UAP are invariant under Banach space isomorphisms.

The partial sum projections of a Schauder basis are uniformly bounded and converge to the identity uniformly on compact sets, giving the BAP. Orthogonal projections give the MAP and the UAP on Hilbert spaces. By Grothendieck's Theorem \cite{G}, the AP, BAP and MAP coincide on reflexive spaces. Even in this class the MAP does not imply the UAP, as Szankowski's uniformly convex example with a symmetric basis shows \cite{Sz76}. The spaces $\ell_p$, $1\le p<\infty$, have the MAP by coordinate truncation and also have the UAP \cite[proof of Theorem~2.1(c)]{FJP}.

In 2009, Szankowski \cite{Sz09} constructed, for each $1<q<2$, a copy of $\ell_q$ in $\ell_q$ whose quotient fails the AP. For $2<p<\infty$ he constructed a subspace of $\ell_p$ which fails the AP and whose quotient is isomorphic to $\ell_p$. In his concluding remark he asked whether the roles of the two exponent ranges could be interchanged. He also asked whether $\ell_p/Y$ must have the BAP whenever $2<p<\infty$ and $Y\subset\ell_p$ has the BAP \cite[p.~282]{Sz09}.

In their paper published in 2011, Figiel, Johnson and Pe{\l}czy\'nski \cite{FJP} asked two questions about approximation in a space with a distinguished closed subspace. The first concerns simultaneous approximation on a space and a subspace. Let $Y$ be a closed subspace of $X$. In the definition of the AP for the pair $(X,Y)$ one requires in addition that $T(Y)\subseteq Y$; a uniform bound on $\|T\|$ gives the BAP for the pair. Such an operator restricts to $Y$ and induces an operator on $X/Y$. In particular, if $(X,Y)$ has the BAP, then $X$, $Y$ and $X/Y$ all have the BAP \cite[Corollary~1.2]{FJP}. Figiel, Johnson and Pe{\l}czy\'nski asked whether the converse holds.
\begin{problem}[{\cite[Problem~6.1]{FJP}}]\label{prob:pair}
Let $Y$ be a closed subspace of $X$. If $X$, $Y$ and $X/Y$ have the BAP, must $(X,Y)$ have the BAP?
\end{problem}
When $Y$ is complemented in $X$, write $X=Y\oplus Z$ with $Z\simeq X/Y$. Taking direct sums of bounded finite-rank approximations to the identities on $Y$ and $Z$ gives uniformly bounded approximations on $X$ that preserve $Y$. Thus the answer is positive in the complemented case.

The authors singled out the case in which $X$ is reflexive. In this case the AP for $(X,Y)$ implies the $1$-BAP for the pair, by Lissitsin and Oja \cite[Corollary~5.12]{LO}. Thus in the reflexive case Problem~\ref{prob:pair} also asks whether the MAP on each of $X$, $Y$ and $X/Y$ implies the AP for the pair.

The second question concerns passage from a subspace to a quotient. Recall that $X$ is an $\mathcal L_p$-space if there is a $\lambda<\infty$ such that each finite-dimensional subspace of $X$ is contained in a finite-dimensional $E\subset X$ with $d(E,\ell_p^{\dim E})\le\lambda$, where $d$ is the Banach--Mazur distance \cite{LP,LR}. Every $\mathcal L_p$-space, $1\le p\le\infty$, has the UAP and hence the BAP and the AP \cite[proof of Theorem~2.1(c)]{FJP}. For $1<p<\infty$ the space $X$ is isomorphic to a complemented subspace of some $L_p(\mu)$ \cite{LR}, and the UAP of $L_p(\mu)$ passes to complemented subspaces. If $X$ is an $\mathcal L_\infty$-space, then the BAP and the UAP each pass from $Y$ to $X/Y$ \cite[Theorem~2.1(a)]{FJP}. For $p=2$ the corresponding assertions follow from Hilbert space geometry and invariance under isomorphisms. This leads to the following question for $2<p<\infty$.
\begin{problem}[{\cite[Problem~6.2]{FJP}}]\label{prob:quotient}
Let $2<p<\infty$, let $X$ be an $\mathcal L_p$-space and let $Y\subset X$ be closed. If $Y$ has the BAP, must $X/Y$ have the BAP? If $Y$ has the UAP, must $X/Y$ have the UAP?
\end{problem}
For $X=\ell_p$, the BAP part of this problem is precisely Szankowski's Question. To see the dual formulation, take $1/p+1/q=1$. Then $X^*$ is an $\mathcal L_q$-space \cite{LR}. In the dual exact sequence the subspace is $Y^\perp=(X/Y)^*$ and the quotient is $X^*/Y^\perp=Y^*$, under the canonical isometries. All three spaces are reflexive. The BAP passes in both directions between a reflexive space and its dual, and the UAP passes in both directions under duality \cite[proof of Theorem~2.1(b)]{FJP}. Thus for $1<q<2$ the equivalent question asks whether the property passes from the quotient to the subspace.

In subsequent work Figiel and Johnson proved that $(X,Y)$ has the BAP when $X$ has the BAP and $Y$ is an $\mathcal L_\infty$-space \cite[Proposition~2.6]{FJ}. When $X$ has the AP, they also characterized the AP for $(X,Y)$ in terms of traces of nuclear operators on $X$ \cite[Theorem~2.2]{FJ}. We use one implication of this characterization in our answer to Problem~\ref{prob:pair}.

A finite-dimensional decomposition of a space is a sequence of finite-dimensional subspaces in which every vector has a unique norm-convergent expansion. In Problem~\ref{prob:pair} the quotient is assumed to have the BAP, and one asks for simultaneous approximation on the pair. In Problem~\ref{prob:quotient} one asks whether the quotient inherits the BAP or the UAP from the subspace. Our first example shows that even finite-dimensional decompositions of both $Y$ and $X/Y$ can coexist with failure of simultaneous approximation in a uniformly convex ambient space.
\begin{theorem}\label{pair}
There is a separable uniformly convex space $X$ with a $1$-unconditional basis and a closed subspace $Y$ such that $Y$ and $X/Y$ have finite-dimensional decompositions. Each of $X,Y,X/Y$ has the MAP, but $(X,Y)$ fails the AP.
\end{theorem}

The next theorem gives the two examples obtained by interchanging Szankowski's exponent ranges.
\begin{theorem}\label{main}
For every $2<p<\infty$ there is an isomorphic embedding $J_p:\ell_p\to\ell_p$ such that $Q_p=\ell_p/J_p(\ell_p)$ fails the AP. For $q=p/(p-1)$, the map $J_p^*:\ell_q\to\ell_q$ is onto and its kernel fails the AP.
\end{theorem}

The range $Y_p=J_p(\ell_p)$ is itself a copy of $\ell_p$, so it has the UAP and hence the BAP. Since $\ell_p/Y_p$ fails the AP, this example answers Szankowski's Question and both parts of Problem~\ref{prob:quotient} in the negative. For each fixed $p$, the same $Q_p$ occurs, up to isomorphism, as a quotient of every infinite-dimensional $\mathcal L_p$-space.

\begin{corollary}\label{allambient}
Let $2<p<\infty$ and let $V$ be an infinite-dimensional $\mathcal L_p$-space. There is an isomorphic embedding $J_V:V\to V$ with range $W=J_V(V)$. The subspace $W$ is isomorphic to $V$ and has the UAP and hence the BAP. The quotient $V/W$ is isomorphic to $Q_p$ and fails the AP.

The adjoint $J_V^*:V^*\to V^*$ is onto and its kernel is isomorphic to $Q_p^*$ and fails the AP.

For $1<q<2$, every infinite-dimensional $\mathcal L_q$-space $V$ admits a surjection $T_V:V\to V$ such that $\ker T_V$ fails the AP and $V/\ker T_V$ has the UAP and hence the BAP.
\end{corollary}
In the range $2<p<\infty$ the proof also gives a closed subspace $W_0\simeq\ell_p$ of $V$ with the UAP for which $V/W_0$ fails the AP. For $V=\ell_p$ one can take $J_V=J_p$, recovering both assertions of Theorem~\ref{main}. We prove the sequence-space case first. The general case follows by applying $J_p$ on a complemented copy of $\ell_p$ in $V$ and the identity on a complementary subspace.

Together with Szankowski's examples, Theorem~\ref{main} shows that $p=2$ is the only exponent in $(1,\infty)$ for which every quotient of $\ell_p$ by a copy of itself has the AP.
\begin{corollary}\label{hilbert}
For $1<p<\infty$ the following are equivalent:
\begin{enumerate}
\item $p=2$.
\item Whenever $Y\subset\ell_p$ is isomorphic to $\ell_p$, the quotient $\ell_p/Y$ has the AP.
\item Every bounded surjection of $\ell_p$ onto itself has a kernel with the AP.
\end{enumerate}
\end{corollary}

The two proofs use the nuclear trace in different ways. For Theorem~\ref{pair} we take an $\ell_2$ sum of finite-dimensional blocks and use random orthogonal identifications to join subspaces of successive blocks. Two subsequences of the block truncations give uniformly bounded approximations on $Y$ and on $X/Y$, respectively. On $X$ we construct a nuclear operator $U$ of trace one with $U(X)\subset Y$ and $U|_Y=0$. For each $T\in\mathcal F(X)$ preserving $Y$, the finite-rank operator $TU$ has square zero, so its ordinary trace is zero. This separates the identity from all such $T$, even though each of the three spaces has the MAP.

For Theorem~\ref{main} we realize $Y$ as the graph of an operator between spaces built from finite-dimensional Hilbert blocks. Estimates for sums of independent random vectors give $X\simeq Y\simeq\ell_p$. On $X/Y$ we construct finite tensors of trace one, indexed by the levels of a binary tree. Cancellation between successive levels and an estimate on the vectors lifting the graph generators give convergence in the projective tensor norm. For each vector in a fixed dense subspace of $X/Y$, the associated operators annihilate that vector at all sufficiently large levels. Thus the limit tensor represents zero and still has trace one.

The paper is organized as follows. In Section~\ref{sec:trace},  we show how the nuclear trace detects failure of the AP for a space and for a pair. In Section~\ref{sec:pair}, we prove Theorem~\ref{pair}. In Section~\ref{sec:graph}, we construct the graph and identify both spaces with $\ell_p$. In Section~\ref{sec:quotient}, we prove convergence of the quotient tensors and complete the proof of Theorem~\ref{main} and Corollaries~\ref{allambient} and~\ref{hilbert}. In Section~\ref{sec:quant}, we bound the ranks of the UAP operators on $Y_p$ and show that both $Q_p$ and $Q_p^*$ have the AP of order $s_p=2p/(3p-2)$.
\section{Trace criteria for failure of the approximation property}\label{sec:trace}
For a single space,  we shall use a tensor of nonzero trace representing the zero operator. For a pair the associated operator can be nonzero: its range must lie in the distinguished subspace, and it must vanish there. Both assertions follow by pairing the tensor with finite-rank approximations to the identity, as in \cite{G,FJ}.

If $A$ and $B$ are Banach spaces, the projective tensor norm of a finite tensor is
\[
 \|u\|_\pi=\inf\left\{\sum_{j=1}^n\|a_j\|\|b_j\|:
                    u=\sum_{j=1}^n a_j\otimes b_j\right\}.
\]
Its completion is denoted by $A\widehat\otimes_\pi B$. Each element admits a series representation with $\sum_j\|a_j\|\|b_j\|<\infty$ \cite{G}. The distinction between a tensor and its associated operator will matter in Section~\ref{sec:quotient}.
For a Banach space $V$ over $\mathbb K$, let $u=\sum_{j=1}^n v_j\otimes\phi_j\in V\otimes V^*$. Using the linear dual over either scalar field, set
\[
 \widetilde u(x)=\sum_{j=1}^n\phi_j(x)v_j,\qquad
 \tr u=\sum_{j=1}^n\phi_j(v_j).
\]
Both $\|\widetilde u\|_{\rm op}$ and $|\tr u|$ are bounded by $\sum_j\|v_j\|\|\phi_j\|$. Taking the infimum over all representations of $u$ gives
\[
 \|\widetilde u\|_{\rm op}\le\|u\|_\pi,\qquad |\tr u|\le\|u\|_\pi .
\]
These bounds extend both maps uniquely by continuity to $V\widehat\otimes_\pi V^*$.

We shall use the following identities for bounded maps $A:V\to W$ and $B:W\to V$ with $BA=I_V$. For $v\otimes\phi$, the tensor $Av\otimes B^*\phi$ has trace $\phi(BAv)=\phi(v)$ and acts on $w\in W$ by $w\mapsto\phi(Bw)Av$. Since $\|Av\|\,\|B^*\phi\|\le\|A\|\|B\|\|v\|\|\phi\|$, the projective norm bound and density give, for every $u\in V\widehat\otimes_\pi V^*$,
\begin{equation}\label{tensortransfer}
 \tr((A\otimes B^*)u)=\tr u,\qquad
 \widetilde{(A\otimes B^*)u}=A\widetilde u B .
\end{equation}

The range of $u\mapsto\widetilde u$ consists of the nuclear operators on $V$. More generally, a nuclear map $L:V\to W$ has a representation $Lx=\sum_j\phi_j(x)w_j$ with $\sum_j\|w_j\|\|\phi_j\|<\infty$. Its nuclear norm $\nu(L)$ is the infimum of these sums. In finite dimensions $\nu(L)=\nu(L^*)$, by transposing a rank-one decomposition in each direction.

For a finite-rank operator $A$, write $\tr_{\rm fin}A$ for the ordinary trace of its restriction to its range. For $u\in V\widehat\otimes_\pi V^*$ put $U=\widetilde u$ and $\Phi_u(T)=\tr((T\otimes I_{V^*})u)$ for bounded $T:V\to V$.

\begin{lemma}\label{tracelemma}
The functional $\Phi_u$ is continuous in the topology of uniform convergence on compact subsets of $V$, and $\Phi_u(T)=\tr_{\rm fin}(TU)$ for every $T\in\mathcal F(V)$. If $U=0$ and $\tr u\ne0$, then $V$ fails the AP.
\end{lemma}
\begin{proof}
Write $u=\sum_j v_j\otimes\phi_j$ with $\|v_j\|=1$ and $\sum_j\|\phi_j\|<\infty$, omitting zero terms. For an infinite sum choose strictly increasing $N_m$ such that $\sum_{j\ge N_m}\|\phi_j\|\le2^{-m}$. Put
\[
 b_j=1+\sum_{m:N_m\le j}1,\qquad
 C_u=\sum_j b_j\|\phi_j\|
     =\sum_j\|\phi_j\|+\sum_m\sum_{j\ge N_m}\|\phi_j\|
     \le\sum_j\|\phi_j\|+1 .
\]
Then $b_j\to\infty$, so $v_j/b_j\to0$. For a finite sum take $b_j=1$. In either case $K=\{0\}\cup\{v_j/b_j\}_j$ is compact and, for every bounded $T$,
\[
 |\Phi_u(T)|=\Big|\sum_j\phi_j(Tv_j)\Big|\le C_u\sup_{x\in K}\|Tx\|.
\]
For $T\in\mathcal F(V)$ choose a basis $(w_i)_{i=1}^r$ of $T(V)$ and write $Tx=\sum_i h_i(x)w_i$. Since $\|v_j\|=1$,
\begin{align*}
 \sum_j|h_i(v_j)\phi_j(w_i)|
 &\le \|h_i\|\,\|w_i\|\sum_j\|\phi_j\| \qquad(1\le i\le r),\\
 \sum_{i=1}^r\sum_j|h_i(v_j)\phi_j(w_i)|
 &\le\Big(\sum_{i=1}^r\|h_i\|\,\|w_i\|\Big)\sum_j\|\phi_j\|<\infty .
\end{align*}
Thus the two sums may be interchanged. The trace of $TU$ on $T(V)$ equals $\tr_{\rm fin}(TU)$ because $TU$ induces the zero operator on $T(V)/TU(V)$. In the basis $(w_i)$ this gives
\begin{equation}\label{finitepairing}
 \begin{aligned}
 \Phi_u(T)&=\sum_j\sum_{i=1}^r h_i(v_j)\phi_j(w_i)
           =\sum_{i=1}^r\sum_j h_i(v_j)\phi_j(w_i)\\
          &=\sum_{i=1}^r h_i(Uw_i)=\tr_{\rm fin}(TU).
 \end{aligned}
\end{equation}
If $U=0$, \eqref{finitepairing} gives $\Phi_u(T)=0$ for every $T\in\mathcal F(V)$. If $V$ had the AP, continuity would imply $\tr u=\Phi_u(I_V)=0$.
\end{proof}

Suppose $V$ has the AP and $\widetilde u=\widetilde v=U$. Since $\widetilde{u-v}=0$, Lemma~\ref{tracelemma} gives $\tr(u-v)=0$. Thus $\tr u=\tr v$. For the resulting nuclear trace $\tr_V$ \cite{G}, we obtain
\[
 \tr_V U=\tr u,\qquad
 \tr_V(TU)=\Phi_u(T)=\tr_{\rm fin}(TU)\quad(T\in\mathcal F(V)).
\]
Taking the infimum of $|\tr_V U|\le\sum_j\|v_j\|\,\|\phi_j\|$ over the nuclear representations of $U$ gives $|\tr_V U|\le\nu(U)$. For finite-rank $L$ on $V$, choose a bounded projection $P$ onto $L(V)$. Then $PL=L$, so $\tr_V L=\tr_V(PL)=\tr_{\rm fin}(PL)=\tr_{\rm fin}L$. Equation~\eqref{finitepairing} itself holds on every Banach space.

For a closed subspace $Y\subset X$, write $\mathcal F_Y(X)=\{T\in\mathcal F(X):T(Y)\subseteq Y\}$.
\begin{corollary}[A trace formula for a pair]\label{pairtrace}
Suppose $Y$ is a closed subspace of $X$, and $u\in X\widehat\otimes_\pi X^*$ satisfies $\widetilde u(X)\subseteq Y$, $\widetilde u|_Y=0$ and $\tr u\ne0$. Then $(X,Y)$ fails the AP.
\end{corollary}
\begin{proof}
Fix $T\in\mathcal F_Y(X)$ and set $U=\widetilde u$, $F=TU(X)$. Since $T$ has finite rank, $F$ is finite dimensional. The range and kernel assumptions give
\[
 F\subseteq T(Y)\subseteq Y\subseteq\ker U,\qquad
 (TU)|_F=0,\qquad \tr_{\rm fin}(TU)=\tr\big((TU)|_F\big)=0 .
\]
By \eqref{finitepairing}, $\Phi_u(T)=0$, whereas $\Phi_u(I_X)=\tr u$. Use the compact set $K$ and constant $C_u$ constructed from $u$ in the proof of Lemma~\ref{tracelemma}. For every $T\in\mathcal F_Y(X)$,
\begin{equation}\label{pairseparation}
 0<|\tr u|=|\Phi_u(I_X-T)|\le C_u\sup_{x\in K}\|x-Tx\|.
\end{equation}
Thus on this fixed $K$ the approximation error is at least $|\tr u|/C_u>0$. For the dual characterization when $X$ has the AP, see \cite{FJ}.
\end{proof}

\section{Failure of the approximation property for a pair in~a~uniformly~convex~space}\label{sec:pair}
Almost every sequence of frames in the construction below has two properties. The induced truncations yield the BAP on $Y$ and on $X/Y$, and the associated block matrix defines a nuclear operator satisfying the hypotheses of Corollary~\ref{pairtrace}.

\subsection{Blocks and graph subspaces}\label{sec:pairblocks}
Work first over $\R$. Let $p=4$, $q=4/3$, and $a=1/4$. On $\R^m$ use the normalized norms
\[
 \ns{x}{s}=m^{1/2-1/s}\|x\|_{\ell_s^m}.
\]
The dimension is always that of the block containing $x$. The Euclidean pairing identifies the dual with exponent $s'$. These norms increase with $s$, by monotonicity of probability-space norms. In particular
$\ns{x}{1}\le\ns{x}{q}\le\eu{x}\le\ns{x}{p}$.
Write $d_k=2^k$, $m_k=3d_k$, $n_k=2d_k=d_{k+1}$ for $k\ge0$. Let $H_k=\R^{m_k}$ with exponent $s_k$, where
\[
 (s_0,s_1,s_2,s_3,s_4,\ldots)=(p,p,q,q,p,\ldots)
\]
is periodic of period four. The consecutive $q$-blocks at $4j+2,4j+3$ will control the truncations on $Y$. The consecutive $p$-blocks at $4j,4j+1$ become $q$-blocks in the dual and give the same estimate on $Y^\perp$, hence on the quotient truncations. Set $X=(\bigoplus_{k\ge0}(H_k,|\cdot|_{[s_k]}))_{\ell_2}$. The coordinate basis is $1$-unconditional. The two families of blocks are scaled copies of finite-dimensional $\ell_4$ and $\ell_{4/3}$ and have a common modulus of uniform convexity. By Day's Theorem on direct sums \cite[Theorem~3]{Day}, $X$ is uniformly convex. In particular $X$ is separable and reflexive. The projections onto its first finitely many blocks are contractive and converge strongly to the identity. Thus $X$ has the MAP.

Let $O_k=(F_k^A\ F_k^B)$ be the coordinate matrices on $\prod_{k\ge0}(O(m_k),\mathrm{Haar})$, with $d_k$ columns in $F_k^A$ and $n_k$ in $F_k^B$. Set
\[
 A_k=\operatorname{ran}F_k^A,\quad B_k=\operatorname{ran}F_k^B,\quad
 W_k=F_{k+1}^A(F_k^B)^\top,\quad W_k^*=W_k^\top .
\]
Write $P_A,P_B$ for Euclidean orthogonal projections onto the indicated spaces. Define
\[
 c_k=\begin{cases}
 1 & s_k=s_{k+1},\\
 n_k^{a/2} & (s_k,s_{k+1})=(p,q),\\
 n_k^{-a/2} & (s_k,s_{k+1})=(q,p).
 \end{cases}
\]
For $b\in B_k$ let $y_k(b)$ have entries $(b,c_kW_kb)$ in blocks $k,k+1$ and zero elsewhere. Set $Y=\overline{\operatorname{span}}\{y_k(b):k\ge0,\ b\in B_k\}$.
Over $\mathbb C$ use the same real frames, the complex coordinate blocks with the same normalized norms, and the closed complex span of the graph generators. Here and below $A_k=F_k^A\mathbb K^{d_k}$ and $B_k=F_k^B\mathbb K^{n_k}$ over the chosen scalar field $\mathbb K$. The same basis, convexity and truncation arguments apply over $\mathbb C$. In both fields we use the linear dual, so the transpose of a real matrix acts on the dual coordinates.
For $x\in Y$ and $f\in Y^\perp\subset X^*$ the respective relations are
\begin{align}
 P_{A_{k+1}}x_{k+1}&=c_kW_kP_{B_k}x_k, \label{primal}\\
 P_{A_{k+1}}f_{k+1}&=-c_k^{-1}W_kP_{B_k}f_k . \label{dual}
\end{align}
Equation~\eqref{primal} holds on the generators of $Y$, hence on $Y$. To prove \eqref{dual}, annihilation of $y_k(b)$ for every $b\in B_k$ gives
\[
 0=\langle b,P_{B_k}f_k+c_k W_k^*P_{A_{k+1}}f_{k+1}\rangle
 \quad(b\in B_k).
\]
The vector paired with $b$ belongs to $B_k$. The coordinate pairing is nondegenerate there over either field, since the real frame satisfies $(F_k^B)^\top F_k^B=I_{n_k}$. Thus the vector vanishes. Applying $W_k$ and using $W_kW_k^*=P_{A_{k+1}}$ gives \eqref{dual}.

\subsection{Uniform truncations for the subspace and quotient}\label{sec:trunc}
The next estimate controls the boundary block in both truncations of Proposition~\ref{individualbap}.

\begin{lemma}\label{matching}
Let $E\subset\R^m$ and $F\subset\R^\ell$ be $n$-dimensional subspaces with $n\le m,\ell\le3n$, and let $W:E\to F$ be a Haar-distributed Euclidean isometry. There exist absolute constants $C,\gamma>0$ such that, with probability at least $1-Ce^{-\gamma n}$, the estimate
\[
 \eu{P_Ex}\le C\big(\ns{x}{1}+\ns{x'}{1}\big)
\]
holds for every $x\in\R^m$ and $x'\in\R^\ell$ satisfying $P_Fx'=\pm WP_Ex$.
\end{lemma}
\begin{proof}
Put $\delta=1/16$. Choose $0<\varepsilon<1$ with $2\varepsilon\log(e(1+6/\varepsilon))\le1/4$, and set $h=\lceil\varepsilon n\rceil$, $t=\delta\sqrt{\varepsilon/3}$. For the Euclidean unit sphere $S_E$, put
\[
 A_E=\{0,\ \pm t\sqrt m\,P_Ee_j:1\le j\le m\},\qquad
 K_E=S_E\cap P_E(t\sqrt m\,B_1^m)=S_E\cap\operatorname{conv}A_E ,
\]
where $B_1^m$ is the $\ell_1^m$ unit ball. For $y\in K_E$ take independent copies $\xi_1,\ldots,\xi_h$ of an $A_E$-valued vector with mean $y$. Then
\[
 \E\left|\frac1h\sum_{i=1}^h\xi_i-y\right|_2^2
 =\frac{\E|\xi_1|_2^2-|y|_2^2}{h}
 \le\frac{t^2m}{h}\le\delta^2 .
\]
The empirical means therefore cover $K_E$ by $\delta$-balls. Choosing one point of $K_E$ from each ball that meets it gives a $2\delta$-net $\mathcal N_E$. Counting multisets gives, for $n\ge1/\varepsilon$,
\[
 |\mathcal N_E|\le\binom{2m+h}{h}
 \le [e(1+2m/h)]^h
 \le [e(1+6/\varepsilon)]^{2\varepsilon n}\le e^{n/4}.
\]
The same construction gives a $2\delta$-net $\mathcal N_F$ of $K_F$ with $|\mathcal N_F|\le e^{n/4}$.

For fixed $v\in S_E$, $w\in S_F$, write $g=g_1w+g_\perp$ for a standard Gaussian in $F$ and put $\widehat g=g/|g|_2$. Then $Wv$ has the law of $\widehat g$. With $r=4\delta=1/4$,
\[
 |\widehat g-w|_2\le r
 \ \Longrightarrow\ g_1>0,\quad |g_\perp|_2\le r|g|_2
 \ \Longrightarrow\ g_1\ge\sqrt{15}\,|g_\perp|_2 .
\]
Independence and the Gaussian tail bound give
\[
 \Pr\{|Wv-w|_2\le r\}
 \le \E e^{-15|g_\perp|_2^2/2}=4^{-(n-1)} .
\]
If $WK_E\cap K_F\ne\varnothing$, some $v\in\mathcal N_E$, $w\in\mathcal N_F$ satisfy $|Wv-w|_2\le4\delta$. Thus
\[
 \Pr\{WK_E\cap K_F\ne\varnothing\}
 \le |\mathcal N_E|\,|\mathcal N_F|\,4^{-(n-1)}
 \le 4e^{-(\log4-1/2)n}.
\]
Enlarging the prefactor covers $n<1/\varepsilon$.

Finally, suppose $P_Fx'=\pm WP_Ex$ and $a=|P_Ex|_2>0$. Since $K_F=-K_F$,
\[
 \max\{\ns{x}{1},\ns{x'}{1}\}\le ta
 \ \Longrightarrow\ W(P_Ex/a)=\pm P_Fx'/a\in WK_E\cap K_F .
\]
Outside the intersection event this gives $a\le t^{-1}(\ns{x}{1}+\ns{x'}{1})$, also valid when $a=0$.
\end{proof}

Fix an edge $k$ and put $\mathcal G_k=\sigma(B_k,A_{k+1})$. Write $\mu_n$ for normalized Haar measure on $O(n)$. For every bounded Borel function $\Psi$ on the space of $m_{k+1}\times m_k$ matrices,
\[
 \E[\Psi(W_k)\mid\mathcal G_k]
 =\int_{O(n_k)}\Psi\bigl(F_{k+1}^A R(F_k^B)^\top\bigr)\,d\mu_{n_k}(R)
 \quad\text{a.s.}
\]
Indeed, the integral depends only on the two subspaces by left and right invariance of $\mu_{n_k}$. For fixed $R$, replacing $O_{k+1}$ by $O_{k+1}\operatorname{diag}(R,I_{n_{k+1}})$ preserves the product measure and both subspaces, since $d_{k+1}=n_k$. Multiplying by any bounded $\mathcal G_k$-measurable function and integrating over $R$ proves the identity. Thus, conditionally on $\mathcal G_k$, the restriction $W_k:B_k\to A_{k+1}$ has the Haar law required in Lemma~\ref{matching}. Put
$\mathcal E_k=\{W_kK_{B_k}\cap K_{A_{k+1}}\ne\varnothing\}$. This event is measurable since its defining constraints are closed and its witnesses range over compact balls. Both sets $K_{B_k},K_{A_{k+1}}$ are fixed under the conditioning. Since $(m_k,m_{k+1})=(3n_k/2,3n_k)$, the uniform estimate in the proof of Lemma~\ref{matching} gives $\Pr(\mathcal E_k\mid\mathcal G_k)\le Ce^{-\gamma n_k}$ almost surely. Hence
\[
 \Pr(\mathcal E_k)=\E\big[\Pr(\mathcal E_k\mid\mathcal G_k)\big]
 \le C e^{-\gamma n_k},\qquad
 \sum_{k\ge0}\Pr(\mathcal E_k)\le C\sum_{k\ge0}e^{-\gamma 2^{k+1}}<\infty .
\]
By the first Borel--Cantelli Lemma, the matching estimate holds for both signs on every sufficiently large edge on the event
\[
 \Omega_{\rm match}=\bigcup_{N\ge0}\bigcap_{k\ge N}\mathcal E_k^c,
 \qquad \Pr(\Omega_{\rm match})=1 .
\]

\begin{proposition}\label{individualbap}
Fix $\omega\in\Omega_{\rm match}$ and choose $k_0$ with $\omega\notin\mathcal E_k$ for every $k\ge k_0$. Let $\rho:X\to X/Y$ be the quotient map and put
\[
 D_kx=(x_0,\ldots,x_k,P_{A_{k+1}}x_{k+1},0,\ldots),\qquad
 \widetilde D_k(\rho x)=\rho D_kx .
\]
These are well-defined finite-rank projections. With $C$ from Lemma~\ref{matching}, set $C_{\mathbb R}=C$, $C_{\mathbb C}=2C$ and $\Lambda=(1+2C_{\mathbb K}^2)^{1/2}$. For $4j\ge k_0$,
\[
 \max\{\|D_{4j+2}|_Y\|,\|\widetilde D_{4j}\|\}\le\Lambda,\qquad
 D_{4j+2}|_Y\longrightarrow I_Y,\quad
 \widetilde D_{4j}\longrightarrow I_{X/Y}
\]
uniformly on compact sets in the respective spaces. In particular, $Y$ and $X/Y$ both have $\Lambda$-BAP.
\end{proposition}
\begin{proof}
Since $P_{A_{k+1}}W_k=W_k$ and $P_{A_{k+1}}B_{k+1}=0$,
\[
 D_k^2=D_k,\qquad D_ky_j(b)=\1_{\{j\le k\}}y_j(b),\qquad D_kY\subseteq Y.
\]
Each $D_k$ has finite block support and is bounded, so it induces the stated quotient projection.

Fix $k\ge k_0$. If $P_{A_{k+1}}v=\sigma W_kP_{B_k}u$, $\sigma=\pm1$, then
\[
 \ns{P_{A_{k+1}}v}{q}\le |P_{B_k}u|_2
 \le C_{\mathbb K}\big(\ns{u}{q}+\ns{v}{q}\big).
\]
Here $|\cdot|_{[1]}\le|\cdot|_{[q]}\le|\cdot|_2$. Over $\mathbb R$, the second inequality follows from Lemma~\ref{matching}. Over $\mathbb C$, the real and imaginary parts satisfy the matching relation separately because the matrices are real, and
$|P_{B_k}u|_2\le |P_{B_k}\Re u|_2+|P_{B_k}\Im u|_2$.

Table~\ref{tab:truncations} applies the estimate to the boundary relations \eqref{primal} and \eqref{dual}. Here $r_i$ denotes the exponent of the $i$th block of $V$.

\begin{table}[htbp]
\caption{Boundary blocks of the primal and adjoint truncations}\label{tab:truncations}
\centering
$\begin{array}{ccccc}\toprule
 k\bmod4 & z\in V & (r_k,r_{k+1},c_k) &
 P_{A_{k+1}}z_{k+1} & T_k z\\ \midrule
 2 & z\in Y\subset X & (q,q,1)& \phantom{-}W_kP_{B_k}z_k & D_k z\\
 0 & z\in Y^\perp\subset X^* & (q,q,1)& -W_kP_{B_k}z_k & D_k^* z\\ \bottomrule
\end{array}$
\end{table}
Since $P_{A_{k+1}}^\top=P_{A_{k+1}}$, the boundary block has the same norm estimate in both cases. Thus
\begin{align*}
 \|T_k z\|_V^2
 &=\sum_{i=0}^k\ns{z_i}{r_i}^2+\ns{P_{A_{k+1}}z_{k+1}}{q}^2\\
 &\le\sum_{i=0}^k\ns{z_i}{r_i}^2
       +2C_{\mathbb K}^2\big(\ns{z_k}{q}^2+\ns{z_{k+1}}{q}^2\big)
 \le\Lambda^2\|z\|_V^2 .
\end{align*}
For $k\ge k_0$ with $k\equiv0\pmod4$, the second row and the onto isometry $\rho^*:(X/Y)^*\to Y^\perp$ give
\[
 \rho^*\widetilde D_k^*=D_k^*\rho^*,\qquad
 \|\widetilde D_k\|=\|D_k^*|_{Y^\perp}\|\le\Lambda .
\]

For each finite sum of graph generators, $D_k$ agrees with the identity for all sufficiently large $k$. The same holds for every $x$ of finite block support. Hence the two sequences converge pointwise to the identity on the respective dense subspaces
\[
 \operatorname{span}\{y_i(b):i\ge0,\ b\in B_i\}\subset Y,\qquad
 \rho(c_{00}(H_i))\subset X/Y .
\]
For either sequence $(S_j)$, take a finite $\delta$-net $F$ for a compact set $K$ with centers in the indicated dense subspace. For all sufficiently large $j$, $S_j|_F=I|_F$ and
\[
 \sup_{z\in K}\|S_jz-z\|
 \le(1+\Lambda)\sup_{z\in K}\min_{v\in F}\|z-v\|
 \le(1+\Lambda)\delta .\qedhere
\]
\end{proof}

\subsection{Nuclear norm estimates for the blocks of \texorpdfstring{$U$}{U}}
We estimate the nuclear norms under the original product Haar measure, then intersect the resulting event of probability one with $\Omega_{\rm match}$ from Proposition~\ref{individualbap}.

For an edge from block $k$ to block $k+1$, set
\[
 h_k x=P_{B_k}x_k-c_k^{-1}W_k^\top x_{k+1},\qquad
 L_k x=n_k^{-1}y_k(h_kx).
\]
Here $h_k:X\to B_k$ vanishes on the generators $y_k(b)$ by the identity $W_k^\top W_k=P_{B_k}$. The range of $L_k$ lies in $y_k(B_k)$. Its only nonzero blocks are on the two coordinates $k,k+1$, where its matrix is
\begin{equation}\label{edgeblock}
 \frac1{n_k}
 \begin{pmatrix}
 P_{B_k}&-c_k^{-1}W_k^\top\\
 c_kW_k&-P_{A_{k+1}}
 \end{pmatrix}.
\end{equation}
We add the matrix entries of the incident edges at each block. For $k\ge1$ the contributions to the diagonal are $-P_{A_k}/n_{k-1}$ and $P_{B_k}/n_k$, with $n_{k-1}=d_k$. At $k=0$ only the second contribution occurs. This gives a tridiagonal matrix with diagonal blocks
\begin{equation}\label{diag}
 V_0=P_{B_0}/n_0,\qquad
 V_k=P_{B_k}/n_k-P_{A_k}/d_k=(I-3P_{A_k})/(2d_k)\quad(k\ge1)
\end{equation}
and two entries on each edge, with rows indexed by output blocks:
\begin{equation}\label{off}
 U_{k+1,k}=c_k W_k/n_k,\qquad
 U_{k,k+1}=-c_k^{-1}W_k^*/n_k .
\end{equation}
For a real matrix $M$ between normalized coordinate blocks, let $\nu_{r,s}(M)$ denote its nuclear norm with source exponent $r$ and target exponent $s$. Under the coordinate pairing,
\[
 \nu_{r,s}(M)=\nu_{s',r'}(M^\top),\qquad
 \nu_{r,s}(v\otimes u)=\ns{v}{s}\ns{u}{r'},
 \quad (v\otimes u)x=(u^\top x)v .
\]

\begin{lemma}\label{diagonalmoment}
For $k\ge1$ and $r\in\{p,q\}$, $\E\nu_{r,r}(V_k)\le C d_k^{-a}$.
\end{lemma}
\begin{proof}
Set $d=d_k$, $m=3d$ and $L=dV_k$. By \eqref{diag},
\[
 L=O_k\operatorname{diag}(-I_d,\tfrac12 I_{2d})O_k^\top,\qquad
 \sum_h\lambda_h=0,\quad \max_h|\lambda_h|=1 ,
\]
where $\lambda_h$ are the diagonal entries in this conjugation. For a Euclidean unit vector $v$, $Q(v)=v^\top Lv$ has the law of $\sum_h\lambda_h u_h^2$ with $u$ uniform on the unit sphere. Write a standard Gaussian vector $g\in\mathbb R^m$ as $g=Ru$, where $R=|g|_2$ is independent of $u$. From $\E(g_h^2-1)^2=2$ and $\E(g_h^2-1)^4=60$,
\begin{align*}
 \E R^8\,\E|Q(v)|^4
 &=\E\Big|\sum_h\lambda_h(g_h^2-1)\Big|^4\\
 &=48\sum_h\lambda_h^4+12\Big(\sum_h\lambda_h^2\Big)^2
 \le48m+12m^2 .
\end{align*}
Since $\E R^8\ge m^4$, polarization at $(e_i\pm e_j)/\sqrt2$ gives
\[
 \max_{i,j}\E|L_{ij}|^4\le C m^{-2},\qquad
 \E\ns{Le_j}{p}^4=m\sum_i\E|L_{ij}|^4\le C .
\]
Using $V_k=d^{-1}\sum_j Le_j\otimes e_j$, $\ns{e_j}{q}=m^{-a}$ and $V_k^\top=V_k$, we obtain
\[
 \E\nu_{q,q}(V_k)=\E\nu_{p,p}(V_k)
 \le \frac{m^{-a}}{d}\sum_{j=1}^m(\E\ns{Le_j}{p}^4)^{1/4}
 \le C\frac{m^{1-a}}{d}\le C' d^{-a} .\qedhere
\]
\end{proof}

\begin{lemma}\label{offest}
For $M\in\{W_k,W_k^\top\}$, $n=n_k$ and $r,s\in\{p,q\}$,
\[
 \E\nu_{r,s}(M/n)\le C b_{r,s}(n),\qquad
 b_{r,s}(n)=
 \begin{cases}
 n^{-a},& r=p\ \text{or}\ s=q,\\
 1,&(r,s)=(q,p).
 \end{cases}
\]
\end{lemma}
\begin{proof}
In either direction write $M=GF^\top$, where $F,G$ are independent Haar $n$-frames in $\mathbb R^m,\mathbb R^\ell$, respectively, with $n\le m,\ell\le3n$. For $u$ uniform on the unit sphere of $\mathbb R^\ell$, the Gaussian radius-direction identity gives
\[
 \E u_i^4=\frac{\E g_i^4}{\E|g|_2^4}=\frac3{\ell(\ell+2)},\qquad
 \E\ns{u}{p}^4=\frac{3\ell}{\ell+2},\qquad \ns{u}{q}\le1 .
\]
Conditionally on $F$, $Me_j$ has the law of $|F^\top e_j|_2\,u$, with $|F^\top e_j|_2\le1$. Hence $\E\ns{Me_j}{s}\le3^{1/4}$ for $s=p,q$. The coordinate expansion and its transpose give
\begin{align*}
 \E\nu_{p,s}(M/n)
 &\le \frac{m^{-a}}{n}\sum_{j=1}^m\E\ns{Me_j}{s}
 \le 3^{1/4}\frac{m^{1-a}}{n}\le C n^{-a},\\
 \E\nu_{r,q}(M/n)&=\E\nu_{p,r'}(M^\top/n)\le C n^{-a}.
\end{align*}
For $(r,s)=(q,p)$, use the columns $u_i$ of $F$ and $v_i$ of $G$. Independence of the two frames and the sphere estimate give
\[
 M=\sum_{i=1}^n v_i\otimes u_i,\qquad
 \E\nu_{q,p}(M/n)
 \le\frac1n\sum_{i=1}^n\E\ns{v_i}{p}\,\E\ns{u_i}{p}\le\sqrt3 .\qedhere
\]
\end{proof}

Write $\nu$ for the nuclear norm with the actual block exponents. For $(r,s)=(s_k,s_{k+1})$ and $n=n_k$, \eqref{off} and Lemma~\ref{offest} bound the two expectations by $C c_k b_{r,s}(n)$ and $C c_k^{-1} b_{s,r}(n)$, respectively. Table~\ref{tab:weights} evaluates both factors for each exponent pair.

\begin{table}[htbp]
\caption{Weights in the nuclear norm estimates for off-diagonal blocks}\label{tab:weights}
\centering
$\begin{array}{cccc}\toprule
 (r,s)&c_k&c_k b_{r,s}(n)&c_k^{-1} b_{s,r}(n)\\ \midrule
 (p,p)\ \text{or}\ (q,q)&1&n^{-a}&n^{-a}\\
 (p,q)&n^{a/2}&n^{-a/2}&n^{-a/2}\\
 (q,p)&n^{-a/2}&n^{-a/2}&n^{-a/2}\\ \bottomrule
\end{array}$
\end{table}
On the product space of frames put
\[
 Z=\sum_{k\ge1}\nu(V_k)+
   \sum_{k\ge0}\big[\nu(U_{k+1,k})+\nu(U_{k,k+1})\big].
\]
The diagonal term in Lemma~\ref{diagonalmoment} uses the marginal law of $O_k$. Both edge terms in Lemma~\ref{offest} use the marginal law $\mu_{m_k}\otimes\mu_{m_{k+1}}$ of $(O_k,O_{k+1})$. Consecutive edges share $O_{k+1}$. Each summand in $Z$ is nonnegative and measurable, so Tonelli's Theorem and the preceding bounds give
\begin{equation}\label{sum}
\begin{aligned}
 \E Z
 &=\sum_{k\ge1}\E\nu(V_k)
   +\sum_{k\ge0}\bigl(\E\nu(U_{k+1,k})+\E\nu(U_{k,k+1})\bigr)\\
 &\le C\sum_{k\ge1}d_k^{-a}+2C\sum_{k\ge0}n_k^{-a/2}\\
 &=C\sum_{k\ge1}2^{-k/4}+2C\sum_{k\ge0}2^{-(k+1)/8}<\infty .
\end{aligned}
\end{equation}
Let $N_{\rm bad}=\sum_{k\ge0}\mathbf 1_{\mathcal E_k}$. By Subsection~\ref{sec:trunc} and \eqref{sum}, $\E(N_{\rm bad}+Z)=\sum_k\Pr(\mathcal E_k)+\E Z<\infty$. Since $\{N_{\rm bad}<\infty\}=\Omega_{\rm match}$, this gives
\begin{equation}\label{jointevent61}
 \Pr\big(\Omega_{\rm match}\cap\{Z<\infty\}\big)
 =\Pr(N_{\rm bad}+Z<\infty)=1 .
\end{equation}
Fix one frame sequence in this intersection. For this sequence, the bounds of Proposition~\ref{individualbap} hold and the block nuclear norms have finite sum. Set $U_{kk}=V_k$ and let $\iota_i,\pi_j$ be the coordinate inclusions and projections. They are contractive, and $\nu(V_0)<\infty$ by finite dimensionality. Consequently
\[
 \sum_{|i-j|\le1}\nu(\iota_i U_{ij}\pi_j)\le\nu(V_0)+Z<\infty,\qquad
 U=\sum_{|i-j|\le1}\iota_i U_{ij}\pi_j
\]
with absolute convergence in nuclear norm. For real coordinate vectors $v,u$, the complex extension of $v\otimes u$ has the same norm product $\ns{v}{s}\ns{u}{r'}$. Taking infima over real block decompositions gives $\nu_{\mathbb C}(U_{ij})\le\nu_{\mathbb R}(U_{ij})$, so the same series converges absolutely over $\mathbb C$.

\subsection{A nuclear operator of trace one}
\begin{proposition}\label{pairblocktrace}
For the spaces $X,Y$ and the nuclear operator $U$ chosen on the event \eqref{jointevent61}, over either scalar field,
\[
 U(X)\subseteq Y,\qquad U|_Y=0,\qquad \tr_X U=1 .
\]
\end{proposition}
\begin{proof}
To determine the action of $U$ on finitely supported vectors we shall sum the finitely many edges meeting their support. To compute its trace we shall instead use the coordinate compressions $S_NUS_N$, where $S_N$ projects onto blocks $0,\ldots,N$. The two sums differ at the last block: compression retains one diagonal part of the boundary edge. Formula~\eqref{boundarytrace} identifies its contribution to the trace.

The real frame matrices act $\mathbb K$-linearly and satisfy $W_k^\top W_k=P_{B_k}$ and $W_kW_k^\top=P_{A_{k+1}}$ over either field. In particular, $h_kx\in B_k$ and $L_k(X)\subseteq y_k(B_k)\subseteq Y$. For $d\in B_j$ the two coordinates of $y_j(d)$ give
\begin{equation}\label{edgeannihilation}
 h_k y_j(d)=
 \begin{cases}
 c_{k-1}P_{B_k}W_{k-1}d=0,&j=k-1,\ k\ge1,\\
 (P_{B_k}-W_k^\top W_k)d=0,&j=k,\\
 -c_k^{-1}W_k^\top d=0,&j=k+1,\\
 0,&|j-k|\ge2 .
 \end{cases}
\end{equation}
The adjacent cases use $W_{k-1}(B_{k-1})=A_k=\ker P_{B_k}$ and $B_{k+1}=\ker W_k^\top$. Thus $L_ky_j(d)=0$ for all $j,k$.

Reading row $r$ from \eqref{edgeblock}, \eqref{diag} and \eqref{off} gives $\pi_0 Ux=h_0x/n_0$ and, for $r\ge1$,
\[
 \pi_r Ux=\frac{c_{r-1}}{n_{r-1}}W_{r-1}h_{r-1}x+\frac{1}{n_r}h_rx .
\]
For $S_Nx$ the last possible output is the $(N+1)$st block of $L_NS_Nx$. Thus
\begin{equation}\label{finiteedgeaction}
 US_N=\sum_{k=0}^N L_kS_N,\qquad
 Uy_j(d)=US_{j+1}y_j(d)=\sum_{k=0}^{j+1}L_ky_j(d)=0 .
\end{equation}
The finite sum $US_Nx$ belongs to $Y$ and $\|Ux-US_Nx\|\le\|U\|\|x-S_Nx\|\to0$. Since $Y$ and $\ker U$ are closed, \eqref{finiteedgeaction} gives
\[
 Ux=\lim_N US_Nx\in Y,\qquad
 Y=\overline{\operatorname{span}}\{y_j(d):j\ge0,\ d\in B_j\}\subseteq\ker U .
\]

For the trace, compress \eqref{finiteedgeaction} on the left. With $P_N=\iota_NP_{B_N}\pi_N$ one has
\[
 S_NL_kS_N=
 \begin{cases}
 L_k,& k<N,\\
 P_N/n_N,&k=N,\\
 0,&k>N .
 \end{cases}
\]
Consequently, for every $N\ge0$ (with the empty sum equal to zero),
\begin{equation}\label{boundarytrace}
 S_NUS_N=\sum_{k<N}L_k+\frac{P_N}{n_N},\qquad
 P_N^2=P_N,\quad P_N X=\iota_N B_N,\quad \operatorname{rank}P_N=n_N .
\end{equation}
For $k<N$, \eqref{edgeannihilation} gives $L_k^2=0$, hence $\tr_{\rm fin}L_k=0$. On edge $N$, the compression retains $P_{B_N}/n_N$ and discards $-P_{A_{N+1}}/n_N$. The two projections have the same rank $n_N=d_{N+1}$ over $\mathbb K$. Thus
\[
 \tr_{\rm fin}(S_NUS_N)
 =\sum_{k<N}\tr_{\rm fin}L_k+\frac{\tr_{\rm fin}P_N}{n_N}
 =\frac{\dim B_N}{n_N}=1 .
\]
Equivalently, summing by diagonal blocks gives $\tr_{\rm fin}V_0=1$ and $\tr_{\rm fin}V_r=n_r/n_r-d_r/d_r=0$ for $r\ge1$. The compression retains the blocks with $i,j\le N$. For the chosen frames, the sum of the nuclear norms of the remaining blocks gives the bound
\begin{equation}\label{compressiontrace}
 \nu(U-S_NUS_N)
 \le\sum_{k>N}\nu(V_k)
   +\sum_{k\ge N}\bigl(\nu(U_{k+1,k})+\nu(U_{k,k+1})\bigr)
 \longrightarrow0 .
\end{equation}
Since $X$ has the MAP, its nuclear trace agrees with the finite-rank trace on compressions and satisfies $|\tr_X L|\le\nu(L)$ by Section~\ref{sec:trace}. Therefore
\[
 |\tr_XU-1|\le\nu(U-S_NUS_N)\longrightarrow0,\qquad
 \nu\Big(U-\sum_{k<N}L_k\Big)\ge
 \Big|\tr_XU-\sum_{k<N}\tr_{\rm fin}L_k\Big|=1 .\qedhere
\]
\end{proof}

\begin{proof}[Proof of Theorem~\ref{pair}]
Fix the frames as in \eqref{jointevent61} and work over either scalar field. Subsection~\ref{sec:pairblocks} gives a separable reflexive space $X$ with the MAP and a closed subspace $Y$. Proposition~\ref{individualbap} gives the BAP on $Y$ and $X/Y$. Let $u\in X\widehat\otimes_\pi X^*$ represent the nuclear operator $U$ of Proposition~\ref{pairblocktrace}. Since $X$ has the MAP, Section~\ref{sec:trace} gives $\tr u=\tr_X U=1$. By $U(X)\subseteq Y$ and $U|_Y=0$, the estimate \eqref{pairseparation} of Corollary~\ref{pairtrace} applies. For its compact set $K\subset X$,
\[
 \inf_{T\in\mathcal F_Y(X)}\ \sup_{x\in K}\|Tx-x\|
 \ge \frac{1}{C_u}>0 .
\]
This rules out the AP for the pair. It remains to check the assertions about $Y$ and $X/Y$ individually. Both are reflexive since $X$ is reflexive. By Grothendieck's Theorem \cite{G}, the BAP on each therefore implies the MAP. Finally $D_iD_j=D_jD_i=D_{\min(i,j)}$. The uniformly bounded subsequences in Proposition~\ref{individualbap} give increasing finite-rank projections converging strongly to the identity. Their successive differences, with the first range included, give the two finite-dimensional decompositions.

On the ambient space one has $\|D_k\|_{X\to X}\to\infty$. Indeed, any uniformly bounded subsequence would fix each finitely supported block vector eventually. By density and the uniform bound it would converge to $I_X$ uniformly on compact subsets of $X$. Since each $D_k$ has finite rank and preserves $Y$, this would give the AP for $(X,Y)$.
\end{proof}

Let $i:Y\hookrightarrow X$ be the inclusion and let $\pi:X\to X/Y$ be the quotient map. Since $U(X)\subseteq Y$ and $U|_Y=0$, the formula $A(\pi x)=Ux$ defines a bounded operator $A:X/Y\to Y$ with $U=iA\pi$. Although $U$ is nuclear on $X$, the map $A\pi:X\to Y$ cannot be nuclear. Indeed, if $A\pi$ had a nuclear representation $\sum_j\phi_j(\,\cdot\,)y_j$, the AP on $X$ and on $Y$ would give
\[
 1=\tr_X U=\sum_j\phi_j(iy_j)=\tr_Y(A\pi i)=0 ,
\]
where the last equality follows from $\pi i=0$.

\section{A graph construction for an embedding of \texorpdfstring{$\ell_p$}{lp} into itself}\label{sec:graph}
We now turn to Problem~\ref{prob:quotient}. Throughout this section and the next, $X$ and $Y$ denote the spaces in a new construction. We shall realize $Y$ as a graph in $X$ and show that each is isomorphic to $\ell_p$. Section~\ref{sec:quotient} will show that $X/Y$ fails the AP.
\subsection{Estimates for finite sums of independent variables}
Fix $2<p<\infty$, let $q=p/(p-1)$ and put $\delta=1/2-1/p$. We use the following consequence of Rosenthal's inequality \cite{R}. For independent mean-zero random vectors $V_i$ in a finite-dimensional Hilbert space,
\begin{equation}\label{ros}
 \mathbb E\Big\|\sum_i V_i\Big\|^p
 \le C_p^p\left(\sum_i\mathbb E\|V_i\|^p+
                   \Big(\sum_i\mathbb E\|V_i\|^2\Big)^{p/2}\right).
\end{equation}
Indeed apply the scalar inequality to the real inner products $\langle g,V_i\rangle$ for fixed $g$, and then integrate over a standard Gaussian $g$ on the underlying real Hilbert space. If $a_p=(\mathbb E|g_1|^p)^{1/p}$, Minkowski's inequality bounds the $L_{p/2}(g)$ norm of $\sum_i\mathbb E_{V_i}|\langle g,V_i\rangle|^2$ by $a_p^2\sum_i\mathbb E\|V_i\|^2$. Dividing by $a_p^p$ gives \eqref{ros}.

The role of $\delta$ can be seen from the estimate in Section~\ref{sec:quotient}. We group the input coordinates into orbits. For an orbit of $N$ inputs in distinct cells, a signed sum of the input vectors has norm $N^{1/p}$ in the ambient $\ell_p$ sum and norm at most $N^{1/p}$ in the quotient. At each fixed nonzero output offset the $N$ coordinates lie in one Euclidean cell, so their signed sum has dual norm at most $N^{1/2}$ times the largest coefficient modulus. Averaging over the signs and summing over the offsets gives an upper bound $C N^{1/p+1/2}=C N^{1-\delta}$ for the projective norm of the orbit tensor, with an absolute constant $C$. The spacing $B$ below lets this gain absorb the number of sample choices. The parameter $m$ controls the terms involving graph lifts.

Set $A_p=2^{1/p}C_p$. Choose a power of two $m$ and a multiple $B$ of eight so large that, with $\eta=m^{-\delta}$,
\begin{equation}\label{choices}
 A_p\sqrt\eta<1,\qquad (1+\delta)/B<\delta/8 .
\end{equation}
Let $d=m^2$, $\mathcal D=\{1,\ldots,m\}\times\{0,\ldots,m-1\}$. Set $L=m(\log_2m+1)$ and give each $\Omega_j=\{0,1\}^j$, $0\le j\le L$, the uniform probability measure, with $\Omega_0$ a singleton. Take $\Omega=\Omega_L$. Read the first $m\log_2m$ bits in successive groups of $\log_2m$ as the binary labels $L_1,\ldots,L_m\in\{0,\ldots,m-1\}$. Read the remaining $m$ bits as signs $\varepsilon_1,\ldots,\varepsilon_m$, with zero denoting $+1$. The labels and signs are mutually independent and each is uniform on its range. Put
\[
 X_{b,l}=m^{1/p}\varepsilon_b\mathbf1_{\{L_b=l\}},\qquad Y_s=\eta^{-2}X_s .
\]
For $\tau\in\Omega_j$ define $X_s^{(j)}(\tau)$ by extending $\tau$ to a word of length $L$ with zeros and evaluating $X_s$ on that word. Set
\[
 K_j(\tau,\tau')=(d\eta^2)^{-1}\sum_s X_s^{(j)}(\tau)X_s^{(j)}(\tau').
\]
Let $\ell_b(\tau)$ and $\sigma_b(\tau)$ be the label and sign of group $b$ read from the zero-filled prefix $\tau$. Summing first over the labels gives
\[
 K_j(\tau,\tau')=\frac1m\sum_{b=1}^m
    \sigma_b(\tau)\sigma_b(\tau')\mathbf1_{\{\ell_b(\tau)=\ell_b(\tau')\}}.
\]
Indeed $m^{2/p}/(d\eta^2)=1/m$. Appending the same label bit preserves equality of the two labels, and appending the same sign bit preserves the product of the two signs. Thus $K_0=1$, all diagonal entries equal one, $|K_j|\le1$, and appending the same bit to both prefixes leaves $K_j$ unchanged. Also $\mathbb E X_sX_r=\eta^2\mathbf1_{\{s=r\}}$, $\|Y_s\|_q=1$, and $\|\sum_s\pm Y_s\|_q\le\|\sum_s\pm Y_s\|_2=\sqrt d/\eta$.

For $\boldsymbol s\in\mathcal D^n$ use the products $X_{\boldsymbol s},Y_{\boldsymbol s}$ on $\Omega^n$. Empty products equal one. For an array of Hilbert vectors write
\[
 A(a)=\Big(\sum\|a_{\boldsymbol s}\|^p\Big)^{1/p},\quad
 B_2(a)=\Big(\sum\|a_{\boldsymbol s}\|^2\Big)^{1/2},\quad
 N_\theta(a)=(A(a)^p+\theta^p B_2(a)^p)^{1/p}.
\]
\begin{lemma}\label{lem:synthesis}
With $\theta=\eta^n$ and $\kappa_n=(A_p\eta^{-1/2})^n$, the map $\mathcal S a=\sum X_{\boldsymbol s}a_{\boldsymbol s}$ satisfies $\|\mathcal S a\|_{L_p(H)}\le\kappa_n N_\theta(a)$.
\end{lemma}
\begin{proof}
For one packet apply \eqref{ros} to $m^{1/p}\varepsilon_b a_{b,L_b}$. The sums of the $p$th and second moments are $A(a)^p$ and $\eta^2 B_2(a)^2$. For $D\subseteq[n]$ put $M_D^p=\sum_{\boldsymbol s_D}(\sum_{\boldsymbol s_{D^c}}\|a_{\boldsymbol s}\|^2)^{p/2}$. Iteration gives
\[
 \|\mathcal S a\|_p^p\le C_p^{pn}\sum_{D\subseteq[n]}\eta^{p(n-|D|)}M_D^p.
\]
For the induction step write $\boldsymbol s=(r,\boldsymbol t)$ and set $a^{(r)}_{\boldsymbol t}=a_{r,\boldsymbol t}$ and $\bar a_{\boldsymbol t}=(a_{r,\boldsymbol t})_{r\in\mathcal D}\in\ell_2(\mathcal D;H)$. Denote synthesis over $h$ packets by $\mathcal S_h$. Conditioning on packets $2,\ldots,n$ and applying the one-packet estimate gives
\[
 \mathbb E\|\mathcal S_n a\|^p\le C_p^p
 \left(\sum_r\mathbb E\|\mathcal S_{n-1}a^{(r)}\|^p
       +\eta^p\mathbb E\|\mathcal S_{n-1}\bar a\|_{\ell_2(\mathcal D;H)}^p\right).
\]
For $D\subseteq\{2,\ldots,n\}$, the first term uses
$\sum_r M_D(a^{(r)})^p=M_{D\cup\{1\}}(a)^p$.
In the second term $M_D(\bar a)=M_D(a)$ because $r$ belongs to the inner Hilbert sum. Its extra factor $\eta^p$ changes $\eta^{p(n-1-|D|)}$ to $\eta^{p(n-|D|)}$. Thus the two terms give all subsets of $[n]=\{1,\ldots,n\}$, starting with the identity for $n=0$.

If $k=n-|D|$, then $M_D\le\min(d^{k\delta}A,B_2)$ by H\"older's inequality. Since $d^\delta=\eta^{-2}$, the geometric mean gives $\eta^kM_D\le(AB_2)^{1/2}\le\eta^{-n/2}N_\theta$. Sum over the $2^n$ subsets.
\end{proof}
\begin{lemma}\label{retract}
For finite-dimensional Hilbert $H$, all $N$ and $0<\theta\le1$, $(H^N,N_\theta)$ is uniformly isomorphic to a complemented subspace of a finite $L_p(\Xi;H)$, with bounds depending only on $p$.
\end{lemma}
\begin{proof}
Take independent symmetric variables $\xi_s$ of magnitude $a=\theta^{-2/(p-2)}$ with probability $a^{-p}$ and zero otherwise. Then $\mathbb E|\xi_s|^p=1$, $\mathbb E\xi_s^2=\theta^2$. Set $Tv=\sum_s\xi_s v_s$, $(Cf)_s=\theta^{-2}\mathbb E(\xi_s f)$. Orthogonality gives $CT=I$, and \eqref{ros} gives $\|T\|\le C_p$. By Bessel's inequality, $\theta B_2(Cf)\le\|f\|_2$. For $h_s\in H^*$, conditioning on the magnitudes and averaging the signs gives
\[
 \mathbb E\Big\|\sum_s\theta^{-2}\xi_s h_s\Big\|^q
 \le\mathbb E\Big(\sum_s\theta^{-4}\xi_s^2\|h_s\|^2\Big)^{q/2}
 \le\sum_s\|h_s\|^q .
\]
Here $q<2$ and $\mathbb E|\theta^{-2}\xi_s|^q=1$. Duality yields $A(Cf)\le\|f\|_p$, so $\|C\|\le2^{1/p}$.
\end{proof}
\subsection{Construction of the operator \texorpdfstring{$J$}{J} and its graph}
The graph has two coordinates. Its first coordinate will carry the functionals used to define the quotient tensor; its second will provide bounded estimates for the terms left over when successive levels are subtracted. We first describe the partition used to norm the first coordinate.
Let $I_t=[2^t,2^{t+1})\cap\mathbb Z$. Number the bits of a node after its leading one from the left. Positions divisible by $B$ are sample positions. Write $h_t=\lfloor t/B\rfloor=nL+j$, $0\le j<L$, and write a node as $(\beta,\gamma)$, where $\beta$ consists of its ordinary bits and $\gamma=(\omega,\tau)\in\Gamma_t=\Omega^n\times\Omega_j$ consists of its sample bits. Thus $\Gamma_t$ has $P=2^{nL+j}=2^{h_t}$ points, each of mass $P^{-1}$. At depth $k\ge2$ put $\ell(v)=v\bmod4$, $r(v)=2\ell(v)+(k\bmod2)$ and partition the nodes into cells with common key
\[
 \mathfrak c(v)=\big(k,\ell(v),(v_i)_{i\in J(k,r(v))}\big),\qquad
 J(k,r)=\{1\le i\le k-2:B\mid i\text{ or }i\equiv r\pmod8\}.
\]
Let $E$ be the $\ell_p$ sum of the Euclidean coordinate spaces on these cells. Write $e_v$ for the coordinate vectors of $E$ and $e_v^*$ for their coordinate functionals. For each $i\in I_t$, $t\ge1$, define the vectors
\[
 a_i=e_{2i}-e_{2i+1},\qquad b_i=\sum_{s=0}^3e_{4i+s},\qquad
 w_i=a_i-b_i/2,\qquad f_i=a_i/2
\]
and the functional
\[
 z_i=e_{2i}^*-e_{2i+1}^*+\sum_{s=0}^3e_{4i+s}^*\in E^* .
\]
The six coefficients of $z_i$ are those of Szankowski's vectors \cite[p.~124]{Sz78}. For $k=i$ the disjoint supports of $a_i$ and $b_i$ give $z_i(w_i)=2-4/2=0$ and $z_i(f_i)=1$. For $k\ne i$ the supports are disjoint unless one node is a child of the other. In that case the two entries of a child difference have opposite signs and are paired with equal coefficients in a four-point sum. Thus $z_i(w_k)=0$ and $z_i(f_k)=\mathbf1_{\{i=k\}}$ for all $i,k$.

At each depth $t$ let $H$ denote the Euclidean space with orthonormal basis $(e_\beta)$ indexed by the ordinary words at that depth, and put
\[
 G_t=L_p(\Gamma_t;H)\oplus_p L_p(\Omega_j;\ell_p^{d^n}(H))
              \oplus_p L_p(\Omega_j;\ell_2^{d^n}(H)).
\]
Take $M_t=\{(r,\kappa_n a,\kappa_n\eta^n a)\}$ in $G_t$, where $r\in L_p(\Gamma_t;H)$ and $a(\tau)=(a_{\boldsymbol s}(\tau))_{\boldsymbol s\in\mathcal D^n}$ is an arbitrary array of vectors in $H$ for each $\tau\in\Omega_j$. Put
\[
 (\mathcal Sa)(\omega,\tau)=\sum_{\boldsymbol s\in\mathcal D^n}
                 X_{\boldsymbol s}(\omega)a_{\boldsymbol s}(\tau).
\]
Apply Lemma~\ref{lem:synthesis} with $\tau$ fixed, take $p$th powers and average over $\tau$. This gives
\[
 \|\mathcal Sa\|_{L_p(\Gamma_t;H)}^p
 \le \kappa_n^p\,\mathbb E_\tau N_{\eta^n}(a(\tau))^p .
\]
Define
\[
 F=r+\mathcal Sa,\qquad \lambda_{\beta,\gamma}=P^{-1/p}F_\beta(\gamma),\qquad
 J_t(r,\kappa_n a,\kappa_n\eta^n a)=\sum_i\lambda_i w_i.
\]
The preceding bound and the triangle inequality give $\|F\|_p\le 2^{1-1/p}\big(\|r\|_p^p+\|\mathcal Sa\|_p^p\big)^{1/p}\le 2^{1-1/p}\|(r,\kappa_n a,\kappa_n\eta^n a)\|_{G_t}$. For each of the six offsets $v=2^h i+s$ in $w_i$, its cell records the sample bits of $i$ through depth $t+h-2$. If $h=1$ and the last bit of $i$ is a sample bit, it is recovered from $\lfloor\ell(v)/2\rfloor$. Thus each output cell for a fixed offset contains inputs with a single sample word. Summing the $p$th powers of the cell norms and using $p/2>1$ bounds the norm at one offset by $(\sum_\gamma(\sum_\beta|\lambda_{\beta,\gamma}|^2)^{p/2})^{1/p}=\|F\|_p$. The sum of the absolute offset coefficients is four. On taking $\ell_p$ sums $M=(\bigoplus_{t\ge1} M_t)_p$, $G=(\bigoplus_{t\ge1} G_t)_p$, the map $Jm=\sum_t J_tm_t$ is bounded, since at each output level only two input levels contribute.

Set $X=E\oplus_pG$, $Y=\{(Jm,m):m\in M\}$, $Q=X/Y$, and denote the quotient map by $\rho$. Each $M_t$ is closed in $G_t$ by finite dimensionality. Since membership in $M$ is determined coordinatewise, $M$ is closed in $G$. Boundedness of $J$ now gives that $m\mapsto(Jm,m)$ is an isomorphism onto the closed graph, with inverse given by the second coordinate. It suffices to identify $X$ and $M$ up to isomorphism.

On a finite-dimensional Hilbert block with orthonormal coordinates, take independent signs on the finite sign cube and define
\[
 Tx=\sum_i\varepsilon_i x_i,\qquad Cf=(\mathbb E\varepsilon_i f)_i.
\]
Orthogonality gives $CT=I$. By \eqref{ros}, $\|Tx\|_p^p\le C_p^p(\sum_i|x_i|^p+\|x\|_2^p)\le 2C_p^p\|x\|_2^p$. Bessel's inequality gives $\|Cf\|_2\le\|f\|_2\le\|f\|_p$. Multiplying each scalar coordinate by the $p$th root of its atom mass identifies the target with a finite-dimensional $\ell_p$ space. Taking $\ell_p$ sums of the maps therefore embeds any countable $\ell_p$ sum of finite Hilbert blocks complementably in $\ell_p$. This applies to $X$.

For $M_t$ the norm to the $p$th power is $\|r\|_p^p+\mathbb E_\tau N_{\eta^n}(\kappa_n a(\tau))^p$. Apply Lemma~\ref{retract} to the coefficient array at each $\tau$, and take the identity on the residual coordinate $r$. The embedding and its left inverse have bounds independent of $t$, so their direct sums embed $M$ complementably into a second copy of $\ell_p$. By Pe{\l}czy\'nski's Theorem \cite{Pel}, every infinite-dimensional complemented subspace of $\ell_p$ is isomorphic to $\ell_p$. Applying this first to $X$ and then to $M$ gives $X\simeq Y\simeq\ell_p$. The complemented embedding of $M$ just constructed is used only to identify the isomorphism class of $Y$.

For later use choose lifts linearly in $\lambda$: for $F=P^{1/p}(\lambda_{\beta,\gamma})_\beta$ set $a_{\boldsymbol s}(\tau)=\eta^{-2n}\mathbb E_\omega X_{\boldsymbol s}F$ and $r=F-\mathcal Sa$. Denote the resulting element of $M_t$ for $\lambda=e_i$ by $g_i$. Then $J_tg_i=w_i$ and $\rho(w_i,0)=-\rho(0,g_i)$. Fix an ordinary word $\beta$. For $y\in L_q(\Gamma_t)$ define a functional on $E$ by
\[
 \zeta_\beta(y)=P^{-1/q}\sum_{\gamma\in\Gamma_t} y(\gamma)z_{\beta,\gamma}.
\]
For any one of the six offsets, distinct sample words give coordinates in distinct cells. The corresponding part of $\zeta_\beta(y)$ therefore has $E^*$ norm
$P^{-1/q}(\sum_\gamma |y(\gamma)|^q)^{1/q}=\|y\|_{L_q(\Gamma_t)}$.
The triangle inequality in $E^*$ gives $\|\zeta_\beta(y)\|_{E^*}\le6\|y\|_{L_q(\Gamma_t)}$.
Since $z_i(w_k)=0$ for all $i,k$, each $(z_i,0)$ annihilates $Y$. We also write $z_i$ for the functional it induces on $Q$. The sum $(\zeta_\beta(y),0)$ thus belongs to $Y^\perp$, and the isometry $\rho^*:Q^*\to Y^\perp$ gives
\begin{equation}\label{dualbound}
 Z_\beta(y)=P^{-1/q}\sum_\gamma y(\gamma)z_{\beta,\gamma},\qquad
 \|Z_\beta(y)\|_{Q^*}=\|\zeta_\beta(y)\|_{E^*}\le6\|y\|_{L_q(\Gamma_t)} .
\end{equation}

\section{The quotient tensor}\label{sec:quotient}
Write $Q^*=Y^\perp$ via the adjoint of $\rho$. We construct a tensor $u$ on $Q$ as a limit of finite tensors $u^{(t)}$. Their traces will all equal one. To pass to the limit we need convergence in the projective norm, so we estimate the differences $u^{(t)}-u^{(t-1)}$ before examining their associated operators.
\subsection{Cancellation between adjacent levels}
For two nodes $i,i'\in I_t$ with different ordinary words put $C_t(i',i)=0$. For equal ordinary words put
\[
 C_t(i',i)=\prod_{\nu=1}^n K_L(\omega_\nu,\omega'_\nu)\ K_j(\tau,\tau').
\]
The coefficient matrices $C_t$ are symmetric, with diagonal entries equal to one and all entries of modulus at most one. For parents $a,a'\in I_{t-1}$ denote their children by $2a+c$ and $2a'+c'$, where $c,c'\in\{0,1\}$. If the position $t$ is a sample bit, appending the same last bit leaves the coefficient unchanged, including when a packet is completed, since $K_0=1$. For an ordinary position $t$ the child entry is $\mathbf1_{\{c=c'\}}C_{t-1}(a',a)$. Set
\[
 Z_{t,i}=\sum_{i'\in I_t} C_t(i',i)z_{i'},\qquad
 u^{(t)}=2^{-t}\sum_{i\in I_t}\rho(f_i,0)\otimes Z_{t,i}
     \in Q\widehat\otimes_\pi Q^* .
\]
By the biorthogonality from Section~\ref{sec:graph} and the diagonal of $C_t$,
\begin{equation}\label{leveltrace}
 \tr u^{(t)}=2^{-t}\sum_{i,i'\in I_t}C_t(i',i)z_{i'}(f_i)
        =2^{-t}\sum_{i\in I_t}C_t(i,i)=1 .
\end{equation}
Since $f_a=b_a/4+w_a/2$, the graph relation from Section~\ref{sec:graph} gives $\rho(f_a,0)=\rho(b_a/4,0)-\rho(0,g_a)/2$. For $t\ge2$ and $v\in I_{t+1}$ put $k(v)=\lfloor v/2\rfloor$, $a(v)=\lfloor v/4\rfloor$ and $\sigma_v=(-1)^v$. Expanding $f_{k(v)}$ in $u^{(t)}$ and $b_{a(v)}$ in $u^{(t-1)}$ yields
\begin{align*}
 u^{(t)}&=2^{-t-1}\sum_{v\in I_{t+1}}\rho(e_v,0)\otimes\sigma_v Z_{t,k(v)},\\
 u^{(t-1)}&=2^{-t-1}\sum_{v\in I_{t+1}}\rho(e_v,0)\otimes Z_{t-1,a(v)}
        -2^{-t}\sum_{a\in I_{t-1}}\rho(0,g_a)\otimes Z_{t-1,a}.
\end{align*}
Subtracting gives
\begin{align}
 u^{(t)}-u^{(t-1)}&=\Delta_t+R_{t-1},\label{diff}\\
 \Delta_t&=2^{-t-1}\sum_{v\in I_{t+1}}\rho(e_v,0)\otimes
                    (\sigma_v Z_{t,k(v)}-Z_{t-1,a(v)}),\nonumber\\
 R_t&=2^{-t-1}\sum_{i\in I_t}\rho(0,g_i)\otimes Z_{t,i}.\nonumber
\end{align}

We estimate the two terms of \eqref{diff} separately. For $\Delta_t$, cancellation within each orbit gives exponential decay in $t$. For the graph lifts in $R_t$, Lemma~\ref{extra} gives decay in the number of complete packets.

\begin{lemma}\label{ordinary}
For every $t\ge2$,
\[
 \|\Delta_t\|_\pi\le 10\,2^{5\delta/4}\,2^{((1+\delta)/B-\delta/8)t}.
\]
In particular, $\sum_{t\ge2}\|\Delta_t\|_\pi<\infty$.
\end{lemma}
\begin{proof}
Write $v=4a+2c+e$, $k=2a+c$ and $\sigma=(-1)^e$. For a mask $\mu$ on the sample positions of $a$ let $a_\mu$ be the word obtained by flipping those bits. Put $p_\mu=C_{t-1}(a_\mu,a)$ and $s_\mu=C_t(2a_\mu+1-c,k)$. Grouping children by their parent and using the prefix identity yields the exact identity
\begin{equation}\label{maskexpand}
 \sigma Z_{t,k}-Z_{t-1,a}
 =\sum_\mu p_\mu\big(\sigma z_{2a_\mu+c}-z_{a_\mu}\big)
 +\1_{\{B\mid t\}}\sum_\mu s_\mu\,\sigma z_{2a_\mu+1-c}.
\end{equation}
The map $\mu\mapsto a_\mu$ enumerates exactly the parents with the same ordinary word as $a$. All other parent coefficients vanish by the definition of $C_{t-1}$, as do their child coefficients in $C_t$. The coefficient $s_\mu$ in the second sum corresponds to the child with final bit $1-c$ and satisfies $|s_\mu|\le1$. For $a'=a_\mu$, write $k'=2a'+c'$ and $\alpha=\1_{\{c'=c\}}$. Apart from its scalar weight the functional is $\sigma z_{k'}-\alpha z_{a'}$. Its coefficient at $4a'+l$ on the middle level is
\[
 B_l=(-1)^{e+(l\bmod2)}\1_{\lfloor l/2\rfloor=c'}-\alpha .
\]
Table~\ref{tab:middle-coefficients} evaluates $B_l$ at the four residues, starting with the input residue.

\begin{table}[htbp]
\caption{Middle-level coefficients for matched and flipped children}\label{tab:middle-coefficients}
\centering
$\begin{array}{lrrrr}
\toprule
\text{Child}&2c+e&2c+1-e&2(1-c)+e&2(1-c)+1-e\\
\midrule
c'=c,\ \alpha=1&0&-2&-1&-1\\
c'=1-c,\ \alpha=0&0&0&1&-1\\
\bottomrule
\end{array}$
\end{table}
Fix $t$, $\mu$ and the last three input bits. Put $r=r(v)$, and take $F$ to be the first $b$ ordinary positions $i\le t-2$ congruent to $r$ modulo $8$, where
\[
 b=\max\{0,\lfloor(t-2)/8\rfloor-h_t\}.
\]
Each residue class contains at least $\lfloor(t-2)/8\rfloor$ positions through $t-2$, of which at most $h_t$ are sample positions. Fix the bits outside $F$ and let the bits in $F$ vary. Each orbit $\mathcal O$ has $M'=2^b$ inputs. Since $F\subseteq J(t+1,r)$, the $M'$ input cell labels $\mathfrak c(v)$ are distinct.

Table~\ref{tab:output-offsets} gives the outputs of $\sigma z_{2a'+c'}-\alpha z_{a'}$ for $a'=a_\mu$. Here $h$ is the output depth and $\zeta=a'\bmod2$.

\begin{table}[htbp]
\caption{Output coordinates, residues and coefficients at the three depths}\label{tab:output-offsets}
\centering
$\begin{array}{ccccc}\toprule
 h&w&\text{Offset}&\ell(w)&\text{Coefficient}\\ \midrule
 t&2a'+d&0\le d<2&2\zeta+d&-\alpha(-1)^d\\
 t+1&4a'+l&0\le l<4&l&B_l\\
 t+2&8a'+4c'+d&0\le d<4&d&\sigma\\ \bottomrule
\end{array}$
\end{table}
On a fixed orbit $\zeta$ is constant, since the final parent bit is frozen and the mask is fixed. In the first and last rows $h\bmod2\ne(t+1)\bmod2$. In the middle row $B_l\ne0$ implies $l\ne2c+e$. Thus, in every row at every nonzero offset,
\[
 r(w)=2\ell(w)+(h\bmod2)\ne 2(2c+e)+((t+1)\bmod2)=r .
\]
When $B\mid t$, the last row has cutoff $h-2=t$ and records the appended sample bit $c'$. This bit is constant on the orbit, including when $c'=1-c$.

Fix a row and offset. Two inputs $v,v'\in\mathcal O$ can differ only at positions in $F$. The common mask preserves their parent differences, and the fixed output suffix introduces no new differences. Thus their outputs differ at exactly the same positions as the inputs. From the definition of the cell key and the choice of $F$,
\[
 F\cap J(h,r(w))\subseteq
 \{i\in F:B\mid i\}\ \cup\ \{i\in F:i\equiv r(w)\pmod8\}
 =\varnothing .
\]
The output depth, residue and recorded bits are therefore identical throughout the orbit. Distinct inputs still give distinct output coordinates, since every position at which they differ is retained in the parent word. Each fixed nonzero offset consequently occupies $M'$ distinct coordinates of a single Hilbert cell. On the flipped branch the middle and last rows give six distinct depth-residue pairs. On the matched branch the sums of the absolute coefficients in the three rows are $2,4,4$, respectively.

Fix one branch of \eqref{maskexpand} and one mask. On each orbit $\mathcal O$, let $\phi_v$ denote the corresponding weighted functional on the right of \eqref{maskexpand}. For independent signs $\delta_v$,
\[
 \sum_{v\in\mathcal O}\rho(e_v,0)\otimes\phi_v
 =\E_\delta\Big[
    \Big(\sum_{v\in\mathcal O}\delta_v\rho(e_v,0)\Big)\otimes
    \Big(\sum_{v\in\mathcal O}\delta_v\phi_v\Big)\Big].
\]
Since the inputs occupy distinct cells,
\[
 \left\|\sum_{v\in\mathcal O}\delta_v\rho(e_v,0)\right\|_Q
 \le\left\|\sum_{v\in\mathcal O}\delta_v e_v\right\|_E=(M')^{1/p}.
\]
For the flipped branch the dual cell norm gives, for every choice of signs,
\[
 \left\|\sum_{v\in\mathcal O}\delta_v\phi_v\right\|_{Q^*}
 =\left(\sum_{d'=1}^{6}
       \left(\sum_{v\in\mathcal O}|s_\mu(v)|^2\right)^{q/2}\right)^{1/q}
 \le 6^{1/q}(M')^{1/2}.
\]
Here $s_\mu(v)$ denotes the coefficient for input $v$. All weighted functionals have their $E^*$ norm under the isometry $Q^*=Y^\perp$. On the matched branch the triangle inequality in $E^*$ over the offsets gives
\[
 \left\|\sum_{v\in\mathcal O}\delta_v\phi_v\right\|_{Q^*}
 \le(2+4+4)\left(\sum_{v\in\mathcal O}|p_\mu(v)|^2\right)^{1/2}
 \le10(M')^{1/2}.
\]
The sign average thus bounds the projective norm of each orbit tensor by $10(M')^{1-\delta}$. There are $2^{t+1}/M'$ orbits per mask and branch, and the number of choices of mask and branch is
\[
 2^{h_{t-1}}(1+\1_{\{B\mid t\}})=2^{h_t}.
\]
Multiplying the orbit bound by $2^{-t-1}$ and summing gives in $Q\widehat\otimes_\pi Q^*$
\[
 \|\Delta_t\|_\pi\le10\,2^{h_t-\delta b},\qquad
 h_t-\delta b\le\big((1+\delta)/B-\delta/8\big)t+5\delta/4 .
\]
The exponent has negative slope by \eqref{choices}.
\end{proof}
\subsection{The terms involving graph lifts}
In \eqref{diff} the terms $R_t$ involve the vectors $g_i$, whose individual norms need not be uniformly bounded. We group these terms using the product representation of the coefficient matrix before taking projective norms.
\begin{lemma}\label{extra}
If $h_t=nL+j$, then $\|R_t\|_\pi\le6\kappa_n\eta^n$. In particular $\sum_t\|R_t\|_\pi<\infty$.
\end{lemma}
\begin{proof}
Let $B_tg=\rho(0,\iota_tg)$, where $\iota_t$ includes $G_t$ into $G$. This contraction is defined on all three coordinate summands, even though an individual summand need not lie in $M_t$. For each full index $l\in\mathcal D$ (a group and a label), put $v_l=X_l^{(j)}$, $v'_l=(d\eta^2)^{-1}X_l^{(j)}$, and $w_l=\|v_l\|_p\|v'_l\|_q$. For each group the probabilities of the zero-filled labels sum to one, so $\sum_l w_l=1$: each label contributes $m^{2/p}(d\eta^2)^{-1}=1/m$ times its probability. For fixed ordinary word $\beta$ the entries of $C_t$ are
\[
 C_t((\beta,\gamma'),(\beta,\gamma))
 =d^{-n}\sum_{\boldsymbol s,l} X_{\boldsymbol s}(\omega)v_l(\tau)
                      Y_{\boldsymbol s}(\omega')v'_l(\tau').
\]
Write $G_\beta(x)=P^{-1/p}\sum_\gamma x(\gamma)g_{\beta,\gamma}$. For $x=X_{\boldsymbol s}v_l$ the lift has zero residual by orthogonality; its coefficient array $h_{\boldsymbol s}$ is supported at $\boldsymbol s$, where its value is $v_l e_\beta$. Thus $G_\beta(x)=a_{\boldsymbol s}^{(p)}+a_{\boldsymbol s}^{(2)}$, where the two terms have respectively the single coordinates $\kappa_n h_{\boldsymbol s}$ and $\kappa_n\eta^n h_{\boldsymbol s}$. Put $b_{\boldsymbol s}=Z_\beta(Y_{\boldsymbol s}v'_l)$. Inserting the kernel expansion in \eqref{diff} gives, with $N_o=2^t/P$,
\[
 R_t=\frac1{2N_o}\sum_{\beta,l}d^{-n}\sum_{\boldsymbol s}
 B_t(a_{\boldsymbol s}^{(p)}+a_{\boldsymbol s}^{(2)})\otimes b_{\boldsymbol s}.
\]
Indeed the two normalizations contribute $P^{1/p+1/q}=P$. For fixed $\beta,l$, $\|a_{\boldsymbol s}^{(2)}\|=\kappa_n\eta^n\|v_l\|_p$ and $\|b_{\boldsymbol s}\|\le6\|v'_l\|_q$ by \eqref{dualbound} and $\|Y_{\boldsymbol s}\|_q=1$. Summing over $\boldsymbol s$ and multiplying by $d^{-n}$ bounds the Hilbert-coordinate term by $6\kappa_n\eta^n w_l$.

For the other coordinate take independent signs $\epsilon_{i,r}$ and put $\epsilon_{\boldsymbol s}=\prod_i\epsilon_{i,s_i}$. Since $\mathbb E(\epsilon_{\boldsymbol s}\epsilon_{\boldsymbol r})=\prod_i\mathbf1_{\{s_i=r_i\}}$, averaging the tensor product of the two signed sums recovers precisely the sum of the elementary tensors. Before $B_t$ the first signed sum has norm $\kappa_n d^{n/p}\|v_l\|_p$. By \eqref{dualbound} the second has norm at most $6\|v'_l\|_q(\sqrt d/\eta)^n$, since the sum factors over packets. After multiplication by $d^{-n}$ the product is at most $6\kappa_n\eta^n w_l$, since $d^{-\delta}=\eta^2$. Summing in $\beta,l$ gives the asserted bound. Each $n$ occurs at most $BL$ times and $\kappa_n\eta^n=(A_p\sqrt\eta)^n$.
\end{proof}
\begin{proof}[Proof of Theorem~\ref{main}]
By \eqref{diff} and the preceding two lemmas the tensors $u^{(t)}$ converge in the projective norm to a tensor $u$ on $Q$. By \eqref{leveltrace}, $\tr u=1$. Each associated operator $\widetilde{u^{(t)}}$ annihilates $\rho(x,g)$ once $t$ exceeds the support levels of a finitely supported $x$, since the functionals $z_i$ on level $t$ use only depths $t+1$ and $t+2$. Such vectors are dense and $\|\widetilde{u^{(t)}}-\widetilde u\|\le\|u^{(t)}-u\|_\pi$, so $\widetilde u=0$. Lemma~\ref{tracelemma} shows that $Q$ fails the AP.

Choose isomorphisms $S:X\to\ell_p$ and $D:\ell_p\to Y$, and set $J_p=SD$ (using the inclusion of $Y$ into $X$). The induced isomorphism $A:Q\to Q_p$ transports $u$ to the tensor
\begin{equation}\label{eq:quotient-witness}
 u_p=(A\otimes(A^{-1})^*)u\in Q_p\widehat\otimes_\pi Q_p^*,
 \qquad \tr u_p=1,\qquad \widetilde u_p=0 .
\end{equation}
The two identities follow from \eqref{tensortransfer} with $B=A^{-1}$. By the Hahn--Banach Theorem, the adjoint $J_p^*$ is onto. Its kernel identifies with $Q_p^*$. Transposition sends $x\otimes f$ to $f\otimes j_{Q_p}x$, where $j_{Q_p}:Q_p\to Q_p^{**}$ is canonical. By continuity the transposed tensor represents the adjoint of the zero operator and retains trace one. Another application of Lemma~\ref{tracelemma} finishes the proof.
\end{proof}

\subsection{Consequences for the ambient space and the exponent}
\begin{proof}[Proof of Corollary~\ref{allambient}]
Every infinite-dimensional $\mathcal L_r$-space contains a complemented copy of $\ell_r$ for $1<r<\infty$ \cite[(11.16), p.~268]{PB}. Such an $\mathcal L_r$-space also has the UAP \cite[proof of Theorem~2.1(c)]{FJP}.

Let $V$ be as in the first assertion and choose coordinates $V\simeq\ell_p\oplus H$ from a complemented copy of $\ell_p$. In these coordinates set $J_V=J_p\oplus I_H$. Its range $W=Y_p\oplus H$ is isomorphic to $V$, so $W$ has the UAP and hence the BAP. The map $(x,h)\mapsto x+Y_p$ induces an isomorphism $V/W\simeq\ell_p/Y_p=Q_p$. Thus $V/W$ fails the AP by Theorem~\ref{main}. The lower bound for $J_V$ and the Hahn--Banach Theorem give surjectivity of $J_V^*$, and $\ker J_V^*=W^\perp\simeq(V/W)^*\simeq Q_p^*$. The latter space fails the AP by Theorem~\ref{main}.

For the additional choice $W_0=Y_p\oplus\{0\}$ we have $W_0\simeq\ell_p$ and $V/W_0\simeq Q_p\oplus H$. The quotient fails the AP because its complemented subspace $Q_p$ fails the AP.

Finally, for an infinite-dimensional $\mathcal L_q$-space with $1<q<2$, write $V\simeq\ell_q\oplus H$ and set $p=q/(q-1)$. The map $T_V=J_p^*\oplus I_H$ is onto, and $\ker T_V\simeq\ker J_p^*$ fails the AP. The induced isomorphism $V/\ker T_V\simeq V$ gives the UAP and hence the BAP on the quotient.
\end{proof}

\begin{proof}[Proof of Corollary~\ref{hilbert}]
Closed subspaces and quotients of a Hilbert space are Hilbert spaces. For $p>2$, Theorem~\ref{main} disproves (2), and Szankowski's quotient map followed by an isomorphism onto $\ell_p$ disproves (3) \cite{Sz09}. For $1<p<2$, use \cite{Sz09} for (2) and Theorem~\ref{main} for (3).
Szankowski uses real scalars. Complexification preserves his subspaces, quotients and isomorphisms: take real and imaginary parts in the quotient sequence and extend the isomorphism and its inverse complex linearly. If a complexified subspace or quotient had the AP, then on any compact subset of its real part $E$ the maps $\Re T|_E$ would give real finite-rank approximants, since $\|\Re(Tx)-x\|\le\|Tx-x\|$ for $x\in E$.
\end{proof}
\section{Rank bounds and approximation properties of order \texorpdfstring{$s$}{s}}\label{sec:quant}
We return to the embedding $J_p$ of Theorem~\ref{main}. For its range $Y_p$ we give an explicit rank bound in the UAP. For the quotient $Q_p$ and its dual $Q_p^*$, we deduce approximation of order $s_p$ from Corollary~10 of Reinov and Latif \cite{RL}.
\subsection{A rank bound for the range of the embedding}
Figiel, Johnson and Schechtman \cite[p.~80]{FJS} recall an argument of Kwapie\'n giving an upper bound of order $(k/\varepsilon)^{Ck}$ for the rank of an operator on $L_p$ of norm at most $1+\varepsilon$ that fixes a prescribed $k$-dimensional subspace. We give a direct estimate on $\ell_p$ with operator norm at most two and record the constants obtained by passing to $Y_p$.
\begin{proposition}
For the subspace $Y_p=J_p(\ell_p)$ the interpolating operators in the UAP may be chosen of rank at most $(1+4k^2)^{2k}+k$ on each $k$-dimensional subspace, with operator norm at most $2\|J_p\|\,\|J_p^{-1}\|$, where $J_p^{-1}:Y_p\to\ell_p$.
\end{proposition}
\begin{proof}
For $1<s<\infty$, let $(x_i)_{i=1}^k$ be an Auerbach basis of a $k$-dimensional $V_s\subset\ell_s$, and use the Hahn--Banach Theorem to extend its biorthogonal functionals $x_i^*$ to $\ell_s$ with norm one. Set $v=\sum_i|x_i|$ and $\delta=(2k^2)^{-1}$. A maximal $\delta$-separated subset $\mathcal A$ of the scalar unit ball is a $\delta$-net. Its open $\delta/2$-balls are disjoint and lie in the ball of radius $1+\delta/2$, so $|\mathcal A|\le(1+2/\delta)^{\dim_{\mathbb R}\mathbb K}\le(1+4k^2)^2$. Choose $\alpha$ into the net with $|\alpha(z)-z|\le\delta$. Partition all coordinates $j$ with $v(j)>0$ by the tuples
\[
 \left(\alpha(x_1(j)/v(j)),\ldots,\alpha(x_k(j)/v(j))\right)\in\mathcal A^k .
\]
For each nonempty class $B$ denote its tuple by $(\alpha_{B,i})_i$ and put $D_0=\operatorname{span}\{v\1_B\}$. Then
\[
 \dim D_0\le|\mathcal A|^k\le(1+4k^2)^{2k},\qquad
 \Big\|x_i-\sum_B\alpha_{B,i}v\1_B\Big\|_s\le\delta\|v\|_s\le\delta k .
\]
Each $v\1_B$ lies in $\ell_s$ even if $B$ is infinite. The vectors $e_B=v\1_B/\|v\1_B\|_s$ are nonnegative and have disjoint supports. Set $e_B^*(x)=\sum_{j\in B}e_B(j)^{s-1}x(j)$, so $\|e_B^*\|=e_B^*(e_B)=1$ over either field. For $P_0x=\sum_B e_B^*(x)e_B$,
\[
 P_0|_{D_0}=I_{D_0},\qquad
 \|P_0x\|_s^s=\sum_B|e_B^*(x)|^s\le\sum_B\|x\1_B\|_s^s\le\|x\|_s^s .
\]
Thus $u_i=(I-P_0)x_i$ satisfies $\|u_i\|_s\le2\operatorname{dist}(x_i,D_0)\le2\delta k$. The correction $Tx=P_0x+\sum_{i=1}^k x_i^*(x)u_i$ gives
\begin{align*}
 Tx_i&=P_0x_i+u_i=x_i,\qquad \|T\|\le1+\sum_i\|u_i\|_s\le1+2\delta k^2=2,\\
 \operatorname{ran}T&\subseteq D_0+\operatorname{span}\{u_1,\ldots,u_k\},\qquad
 \operatorname{rank}T\le\dim D_0+k .
\end{align*}
Conjugation by $J_p:\ell_p\to Y_p$ gives the assertion with norm bound $2\|J_p\|\|J_p^{-1}\|$ (the inverse is taken on the range).
\end{proof}
\subsection{The AP of order \texorpdfstring{$s$}{s} and complexification}
\begin{corollary}\label{cor:orderp}
Let $2<p<\infty$ and set
\[
 s_p=\left(1+\left|\frac1p-\frac12\right|\right)^{-1}
     =\frac{2p}{3p-2}.
\]
Both $Q_p$ and $Q_p^*$ have $\mathrm{AP}_{s_p}$ and fail the classical AP. In every projective representation of the tensor $u_p$ in \eqref{eq:quotient-witness}, or of its transpose, the sum of the $s_p$-powers of the products of the factor norms diverges.
\end{corollary}

For complex scalars the assertion about $\mathrm{AP}_{s_p}$ follows from Corollary~10 of Reinov and Latif \cite{RL}. To prove the assertion about tensor representations and obtain the same result for real scalars, we recall the definition of $\mathrm{AP}_s$ and describe its behavior under complexification.

Let $0<s\le1$. For a real or complex Banach space $V$, take all duals and tensor products over its scalar field. On a finite tensor $v$ put $\|v\|_s^s=\inf\sum_i(\|\phi_i\|\,\|x_i\|)^s$, where the infimum runs over all representations $v=\sum_i\phi_i\otimes x_i$. Since $\|\cdot\|_\pi\le\|\cdot\|_s$, the identity extends to $\iota_s$ in the factorization
\[
 V^*\widehat\otimes_s V \xrightarrow{\ \iota_s\ }
 V^*\widehat\otimes_\pi V \xrightarrow{\ \tau\ }
 V\widehat\otimes_\pi V^*,\qquad j_s^V=\tau\iota_s .
\]
Here $\tau(\phi\otimes x)=x\otimes\phi$ is the isometric flip. The associated operator is $A_s v=\widetilde{j_s^V v}$, since both maps send $\phi\otimes x$ to $y\mapsto\phi(y)x$. By Definition~2 of Reinov and Latif \cite{RL}, $V\in\mathrm{AP}_s$ means that $A_s v=0$ implies $\iota_s v=0$ in the projective tensor product. Indeed, on such a space,
\[
 w=j_s^V v,\quad \widetilde w=0
 \ \Longrightarrow\ A_s v=0
 \ \Longrightarrow\ \iota_s v=0
 \ \Longrightarrow\ w=\tau\iota_s v=0 .
\]
Since the tensor trace vanishes at zero, for every $w\in V\widehat\otimes_\pi V^*$ we obtain
\begin{equation}\label{ordersupport}
 V\in\mathrm{AP}_s,\quad \widetilde w=0,\quad \tr w\ne0
 \quad\Longrightarrow\quad w\notin\operatorname{ran}j_s^V .
\end{equation}

If a projective representation $w=\sum_i x_i\otimes\phi_i$ satisfies $\sum_i(\|x_i\|\,\|\phi_i\|)^s<\infty$, its flipped partial sums $v_N=\sum_{i\le N}\phi_i\otimes x_i$ satisfy
\[
 \sup_{M>N}\|v_M-v_N\|_s^s\le\sum_{i>N}(\|x_i\|\,\|\phi_i\|)^s\longrightarrow0 .
\]
Thus $v=\lim_N v_N$ exists in $V^*\widehat\otimes_s V$ and $j_s^V v=w$. Consequently, if $w\notin\operatorname{ran}j_s^V$, every projective representation of $w$ has divergent $s$-power sum.

\begin{lemma}[Passage to real scalars]\label{lem:realdescent}
Let $0<s\le1$ and let $E=V+iV$ be a complexification of a real Banach space $V$ such that the inclusion of $V$ and conjugation on $E$ are isometries. For $w\in V\widehat\otimes_{\pi,\mathbb R}V^*$ and $x,y\in V$, the canonical complexification maps $c_s,c_\pi$ satisfy
\begin{equation}\label{complextrace}
 \begin{aligned}
 j_s^E c_s&=c_\pi j_s^V,\qquad \tr(c_\pi w)=\tr w,\\
 \widetilde{c_\pi w}(x+iy)&=\widetilde w(x)+i\widetilde w(y).
 \end{aligned}
\end{equation}
Moreover, $c_\pi$ admits a real-linear left inverse $R$ with $\|R\|\le2$, and
\begin{equation}\label{realdescent}
 E\in\mathrm{AP}_s\ \Longrightarrow\ V\in\mathrm{AP}_s .
\end{equation}
In particular, for a closed subspace $Y\subset\ell_p(\mathbb R)$, $1<p<\infty$, put $Y_c=Y+iY$ and $q=p/(p-1)$. The complex quotient $\ell_p(\mathbb C)/Y_c$ and the complex annihilator $Y_c^\perp\subset\ell_q(\mathbb C)$ satisfy the isometry hypotheses for $V=\ell_p(\mathbb R)/Y$ and $V=Y^\perp\subset\ell_q(\mathbb R)$, respectively.
\end{lemma}
\begin{proof}
Since conjugation on $E$ is isometric, the real and imaginary parts are contractions onto $V$. For $\phi\in V^*$, its complex-linear extension satisfies
\[
 \phi_{\mathbb C}(x+iy)=\phi(x)+i\phi(y),\qquad
 |\phi_{\mathbb C}(z)|=\sup_\theta\phi(\Re(e^{i\theta}z))\le\|\phi\|\,\|z\|.
\]
The real-linear maps
\[
 c_s(\phi\otimes_{\mathbb R}x)=\phi_{\mathbb C}\otimes_{\mathbb C}x,
 \qquad
 c_\pi(x\otimes_{\mathbb R}\phi)=x\otimes_{\mathbb C}\phi_{\mathbb C}
\]
extend to contractions on the respective completions. Indeed, $\|\phi_{\mathbb C}\|=\|\phi\|$ by restriction to $V$ and the preceding estimate, so each finite real representation gives a complex representation with the same sum of $s$-powers of norm products (with exponent $1$ for $c_\pi$). Both $j_s^E c_s$ and $c_\pi j_s^V$ send $\phi\otimes_{\mathbb R}x$ to $x\otimes_{\mathbb C}\phi_{\mathbb C}$. This tensor has trace $\phi(x)$ and sends $a+ib$ to $(\phi(a)+i\phi(b))x$. Density and continuity of the tensor trace and the operator map give \eqref{complextrace}.

To construct a left inverse, write $z=x+iy\in E$ and $F|_V=f+ig$ for $F\in E^*$. The real-bilinear map $\beta(z,F)=x\otimes_{\mathbb R}f-y\otimes_{\mathbb R}g$ satisfies
\[
 \beta(iz,F)=-y\otimes_{\mathbb R}f-x\otimes_{\mathbb R}g=\beta(z,iF).
\]
As a real vector space, $E\otimes_{\mathbb C}E^*$ is the quotient of the real algebraic tensor product of the underlying real spaces by
\[
 N=\operatorname{span}_{\mathbb R}
 \{\,iz\otimes_{\mathbb R}F-z\otimes_{\mathbb R}iF:z\in E,\ F\in E^*\,\}.
\]
The real linearization of $\beta$ vanishes on $N$, so $R_0(z\otimes_{\mathbb C}F)=\beta(z,F)$ is independent of the representation. For every finite representation $t=\sum_k z_k\otimes_{\mathbb C}F_k$, write $z_k=x_k+iy_k$ and $F_k|_V=f_k+ig_k$. Then
\[
 \|R_0t\|_{\pi,\mathbb R}
 \le\sum_k(\|x_k\|\,\|f_k\|+\|y_k\|\,\|g_k\|)
 \le 2\sum_k\|z_k\|\,\|F_k\|.
\]
Taking the infimum over all complex representations gives $\|R_0t\|_{\pi,\mathbb R}\le2\|t\|_{\pi,\mathbb C}$. Thus $R_0$ extends uniquely to a real-linear map
\[
 R:E\widehat\otimes_{\pi,\mathbb C}E^*\longrightarrow
 V\widehat\otimes_{\pi,\mathbb R}V^*,\qquad \|R\|\le2 .
\]
Since $\phi_{\mathbb C}|_V=\phi$, we have
$Rc_\pi(x\otimes_{\mathbb R}\phi)=\beta(x,\phi_{\mathbb C})=x\otimes_{\mathbb R}\phi$.
Density gives $Rc_\pi=I$ on $V\widehat\otimes_{\pi,\mathbb R}V^*$.

Suppose $E\in\mathrm{AP}_s$. Take $v\in V^*\widehat\otimes_{s,\mathbb R}V$ with $w=j_s^V v$ and $\widetilde w=0$. By \eqref{complextrace},
$\widetilde{j_s^E c_s v}=\widetilde{c_\pi w}=0$. Definition~2 of \cite{RL} on $E$ gives $j_s^E c_s v=0$, whence
\[
 j_s^V v=Rc_\pi w=Rj_s^E c_s v=0 .
\]
This proves \eqref{realdescent} by the defining implication on $V$.

For the quotient, $Y_c$ is closed by continuity of real and imaginary parts, and $x+iy\in Y_c$ exactly when $x,y\in Y$. Thus $(x+Y)+i(y+Y)\mapsto(x+iy)+Y_c$ is a complex-linear bijection. For real $x$ and arbitrary complex $z$,
\[
 \begin{aligned}
 \|x+Y_c\|&=\inf_{u,v\in Y}\|(x-u)-iv\|_p=\inf_{u\in Y}\|x-u\|_p=\|x+Y\|,\\
 \|\bar z+Y_c\|&=\inf_{h\in Y_c}\|\bar z-\bar h\|_p=\|z+Y_c\|.
 \end{aligned}
\]
The first equality of norms follows from $\|x-u\|_p\le\|(x-u)-iv\|_p$ with equality at $v=0$. The conjugation identity follows from $\{\bar h:h\in Y_c\}=Y_c$.

For the annihilator, use the complex-linear dual pairing $f(z)=\sum_j f_jz_j$. Since $Y$ spans $Y_c$ over $\mathbb C$, for $a,b\in\ell_q(\mathbb R)$,
\[
 a+ib\in Y_c^\perp
 \ \Longleftrightarrow\ a(y)+ib(y)=0\quad(y\in Y)
 \ \Longleftrightarrow\ a,b\in Y^\perp .
\]
Thus $Y_c^\perp=Y^\perp+iY^\perp$. In its inherited $\ell_q(\mathbb C)$ norm,
\[
 \|a\|_{\ell_q(\mathbb C)}=\|a\|_{\ell_q(\mathbb R)},\qquad
 \|a-ib\|_q^q=\sum_j(a_j^2+b_j^2)^{q/2}=\|a+ib\|_q^q .
\]
\end{proof}

\begin{proof}[Proof of Corollary~\ref{cor:orderp}]
Write $Y_p=J_p(\ell_p)$ and $q=p/(p-1)$. The quotient $Q_p$ is a quotient of $\ell_p$, and $Q_p^*=Y_p^\perp$ is isometric to a subspace of $\ell_q$. Since $|1/q-1/2|=|1/p-1/2|$, Corollary~10 of Reinov and Latif \cite[p.~193]{RL} gives the same order $s_p$ for these two spaces over $\mathbb C$. Over $\mathbb R$ apply that corollary to the complex quotient by $Y_p+iY_p$ and to its annihilator, then apply Lemma~\ref{lem:realdescent} to the two real spaces.

Theorem~\ref{main} gives the failure of the AP for $Q_p$ and for $Q_p^*\simeq\ker J_p^*$. The tensor $u_p$ and its transpose have trace one and represent zero. By \eqref{ordersupport} neither lies in the range of the corresponding $j_{s_p}$, which gives the assertion about representations.
\end{proof}
For $p=4$ the order is $4/5$. As $p$ decreases to $2$, the orders $s_p$ tend to $1$. Thus within this family the guaranteed order of approximation can be arbitrarily close to $1$ even though each quotient fails the AP.

\section*{AI Declaration}
During the preparation of this work, the authors used GPT-6 Astra to assist with language editing and the refinement of mathematical arguments. The authors take full responsibility for the content of this paper.


\begin{thebibliography}{99}
\enlargethispage{2\baselineskip}
\bibitem{Day} M. M. Day, Some more uniformly convex spaces, \emph{Bull. Amer. Math. Soc.} \textbf{47} (1941), 504--507.
\bibitem{FJ} T. Figiel and W. B. Johnson, The dual form of the approximation property for a Banach space and a subspace,
\emph{Studia Math.} \textbf{231} (2015), 287--292.
\bibitem{FJP} T. Figiel, W. B. Johnson and A. Pe{\l}czy\'nski, Some approximation properties of Banach spaces and Banach lattices,
\emph{Israel J. Math.} \textbf{183} (2011), 199--231.
\bibitem{FJS} T. Figiel, W. B. Johnson and G. Schechtman, Factorizations of natural embeddings of $\ell_p^n$ into $L_r$, I,
\emph{Studia Math.} \textbf{89} (1988), 79--103.
\bibitem{G} A. Grothendieck, \emph{Produits tensoriels topologiques et espaces nucl\'eaires},
Mem. Amer. Math. Soc., no.~16, 1955, 196+140~pp.
\bibitem{LO} A. Lissitsin and E. Oja, The convex approximation property of Banach spaces,
\emph{J. Math. Anal. Appl.} \textbf{379} (2011), 616--626.
\bibitem{RL} O. Reinov and Q. Latif, Grothendieck--Lidski\u{\i} theorem for subspaces of quotients of $L_p$-spaces,
\emph{Banach Center Publ.} \textbf{102} (2014), 189--195.
\bibitem{R} H. P. Rosenthal, On the subspaces of $L^p$ $(p>2)$ spanned by sequences of independent random variables,
\emph{Israel J. Math.} \textbf{8} (1970), 273--303.
\bibitem{Sz76} A. Szankowski, A Banach lattice without the approximation property,
\emph{Israel J. Math.} \textbf{24} (1976), 329--337.
\bibitem{Sz78} A. Szankowski, Subspaces without the approximation property,
\emph{Israel J. Math.} \textbf{30} (1978), 123--129.
\bibitem{Sz09} A. Szankowski, Three-space problems for the approximation property,
\emph{J. Eur. Math. Soc.} \textbf{11} (2009), 273--282.
\bibitem{LP} J. Lindenstrauss and A. Pe{\l}czy\'nski, Absolutely summing operators in $\mathcal L_p$-spaces and their applications, \emph{Studia Math.} \textbf{29} (1968), 275--326.
\bibitem{LR} J. Lindenstrauss and H. P. Rosenthal, The $\mathcal L_p$ spaces, \emph{Israel J. Math.} \textbf{7} (1969), 325--349.
\bibitem{PB} A. Pe{\l}czy\'nski (in collaboration with Cz. Bessaga), Some aspects of the present theory of Banach spaces, in S. Banach, \emph{\OE uvres}, vol.~II, PWN, Warsaw, 1979.
\bibitem{Pel} A. Pe{\l}czy\'nski, Projections in certain Banach spaces, \emph{Studia Math.} \textbf{19} (1960), 209--228.
\end{thebibliography}
\end{document}